\documentclass[12pt,a4paper]{article}
\usepackage[T1]{fontenc}
\usepackage{lmodern}
\usepackage[a4paper,margin=2.5cm]{geometry}
\usepackage{amsmath,amssymb,amsthm,mathtools}
\usepackage{microtype,booktabs,array,longtable,capt-of}
\usepackage{graphicx,adjustbox}
\usepackage{enumitem,etoolbox}
\usepackage[numbers,sort&compress]{natbib}
\usepackage{xurl}
\usepackage[hidelinks]{hyperref}
\allowdisplaybreaks[1]
\setlist{itemsep=2pt,topsep=4pt}
\numberwithin{equation}{section}
\newtheorem{theorem}{Theorem}[section]
\newtheorem{lemma}[theorem]{Lemma}
\newtheorem{proposition}[theorem]{Proposition}
\newtheorem{corollary}[theorem]{Corollary}
\theoremstyle{definition}

\newtheorem{problem}[theorem]{Problem}
\newtheorem{remark}[theorem]{Remark}
\newcommand{\R}{\mathbb R}
\newcommand{\one}{\mathbf 1}
\newcommand{\tr}{\operatorname{tr}}
\newcommand{\rank}{\operatorname{rank}}
\newcommand{\diag}{\operatorname{diag}}
\newcommand{\Spec}{\operatorname{Spec}}
\newcommand{\dist}{\operatorname{dist}}
\newcommand{\sgn}{\operatorname{sgn}}

\newcommand{\Frob}[1]{\left\|#1\right\|_F}
\newcommand{\doi}[1]{\href{https://doi.org/#1}{doi:\nolinkurl{#1}}}
\newcommand{\bigsq}{\mathbin{\dot\cup}}
\title{Extremal Ky Fan $p$-$k$ norms of graphs\thanks{Supported by the National Natural Science Foundation of China (No.\ 12661070).}}
\author{Zhen Lin\thanks{Corresponding author. E-mail: \texttt{linzhen@qhnu.edu.cn}.}\\[3pt]
\small School of Mathematics and Statistics, Qinghai Normal University,\\
\small Xining 810008, Qinghai, China}
\date{}
\begin{document}
\maketitle
\begin{abstract}
We study a two-parameter extension of Nikiforov's extremal problems for graph norms, focusing on the quantity
$$
S_{p,k}(G) = \left(\sum_{i=1}^{k}\sigma_i(G)^p\right)^{1/p},
$$
where $p\ge1$, and $\sigma_1(G) \ge \cdots \ge \sigma_n(G)$ are the singular values of the adjacency matrix of an $n$-vertex graph $G$. This quantity provides a common generalization of the Ky Fan $k$-norm and the Schatten $p$-norm, which correspond to the cases $p=1$ and $k=n$, respectively. We establish upper and lower bounds and characterize extremal and near-extremal graphs within specified parameter ranges. For each fixed $k \ge 4$, we determine an explicit exponent $\vartheta_k \in [1,2)$ such that the balanced complete $k$-partite graph is the unique maximizer for all sufficiently large orders, and we characterize the unique graph attaining the second-largest value. In particular, this determines the maximum Ky Fan $4$-norm for all sufficiently large orders. At these exponents, any maximizer with a fixed edge deficit from the Tur\'{a}n graph's edge count belongs to an explicit finite family of complete $k$-partite graphs (whenever such a realization exists), whereas a fixed edge surplus is uniquely maximized by adding a star inside a largest part of the balanced graph. Furthermore, we quantify the relation between small absolute norm deficits and graph structure, demonstrating that relative near-optimality alone does not guarantee structural proximity. We also obtain a vertex-partition lower bound featuring an optimal universal remainder coefficient and complete equality conditions. These results provide structural answers to Nikiforov's norm questions in the two-parameter setting, linking extremal values to uniqueness, prescribed edge counts, and the structure of near-extremal graphs.
\end{abstract}
\noindent\textbf{MSC} 05C50, 05C35, 15A18, 15A60.\\
\textbf{Keywords:} Ky Fan $p$-$k$ norms; Ky Fan norm; Schatten norm; Extremal graph; Singular value.

\section{Introduction}\label{sec:introduction}
Let $G$ be a simple graph with $n$ vertices and $m$ edges. Write
$\lambda_1(G)\ge\cdots\ge\lambda_n(G)$ for the eigenvalues of its adjacency
matrix $A=A(G)$, and $\sigma_1(G)\ge\cdots\ge\sigma_n(G)\ge0$ for its
singular values. Thus the singular values are the absolute values of the
eigenvalues, counted with multiplicity and reordered. For $p\ge1$ and
$1\le k\le n$, put
\begin{equation}\label{eq:def}
 S_{p,k}(G)=\left(\sum_{i=1}^k\sigma_i(G)^p\right)^{1/p}.
\end{equation}
This is the Ky Fan $k$-norm when $p=1$ and the Schatten $p$-norm when
$k=n$. In particular, $S_{1,n}(G)$ is the energy of $G$
\cite{Gutman1978,So2010}.

In \cite[Problem~1]{Nik2011}, Nikiforov proposed studying the extrema
of the Ky Fan norms of graphs and their relation to graph structure.
The corresponding question for classes of matrices appears as
\cite[Problem~2]{Nik2011}. His subsequent work
\cite{Nik2012,Nik2015,Nik2016} treats Ky Fan and Schatten norms within a
common framework. Motivated by these problems, we study the extremal quantities
\[
 M_{p,k}(n)=\max_{|V(G)|=n}S_{p,k}(G),\qquad
 M_{p,k}(n,m)=\max_{\substack{|V(G)|=n\\|E(G)|=m}}S_{p,k}(G).
\]
For $p=1$, $M_{1,k}(n)$ is Nikiforov's function $\xi_k(n)$.
Our aim is to obtain bounds for these quantities and to determine the
extremal graphs in the parameter ranges stated below. 

Mohar~\cite{Mohar2009} studied the sum of the $k$ largest adjacency
eigenvalues in algebraic order. Nikiforov~\cite{Nik2011} obtained the
corresponding singular-value bound
\[
 S_{1,k}(G)\le\frac{1+\sqrt{k}}2\,n.
\]
The distinction between the two sums matters: singular values do not record
the signs of eigenvalues. For $k=2$, Nikiforov's observation linking the
maximum Ky Fan norm to the maximum spread, together with the theorem of
Breen, Riasanovsky, Tait and Urschel~\cite{Breen2022}, gives
\[
 \lim_{n\to\infty}\frac{M_{1,2}(n)}n
 =\frac2{\sqrt3}<\frac{1+\sqrt2}{2}.
\]
Here the relation concerns the two maximization problems. In particular,
the order bound is not asymptotically attained for every truncation index.

For the full trace norm, Koolen and Moulton~\cite{KoolenMoulton2001}
obtained a sharp order bound and characterized equality. For
$1\le p<2$, Nikiforov~\cite{Nik2012,Nik2016} developed the corresponding
second-moment method for Schatten norms. In particular,
\cite[Corollary~4.23 and Proposition~4.24]{Nik2016} give bounds in terms
of the spectral radius and the order and size. The truncated spectral-radius
bound used here is
\begin{equation}\label{eq:intro-perron}
 S_{p,k}(G)^p\le\rho^p+(k-1)^{1-p/2}(2m-\rho^2)^{p/2},
 \qquad \rho=\lambda_1(G)=\sigma_1(G).
\end{equation}
Together with $2m\le n\rho$, it gives $S_{p,k}(G)\le nC_{p,k}$, where
\begin{equation}\label{eq:intro-envelope}
 C_{p,k}^{\,p}=\max_{0\le x\le1}
 \left\{x^p+(k-1)^{1-p/2}[x(1-x)]^{p/2}\right\}.
\end{equation}
At $k=n$, Theorem~\ref{thm:fixed-size} recovers these full-spectrum
bounds; for $k<n$, its equality condition also requires the omitted
singular values to vanish. At $p=1$, the order bound is precisely
Nikiforov's Ky Fan bound, whereas for $1<p<2$ its coefficient strictly
improves the direct centered-matrix estimate
(Proposition~\ref{prop:comparison}). The second-moment argument is classical;
we use its quantitative deficit to recover the adjacency structure of $G$.

For $k\ge4$, define
\begin{equation}\label{eq:critical-exponent}
 \vartheta_k=\frac{\log((k-1)(k-2)/2)}{\log(k-1)},\qquad
 C_k=\frac1k\binom{k}{2}^{1/\vartheta_k}.
\end{equation}
Then $\vartheta_4=1$, $1<\vartheta_k<2$ for $k\ge5$, and $\vartheta_k\to2$. At $p=\vartheta_k$ the
maximum in~\eqref{eq:intro-envelope} occurs at $(k-1)/k$ and equals
$C_k^{\vartheta_k}$. Let $T_k(n)$ denote the complete $k$-partite graph whose part
sizes differ by at most one.

\begin{theorem}\label{main:order}
Fix $k\ge4$. Every graph $G$ of order $n\ge k$ satisfies
$S_{\vartheta_k,k}(G)\le nC_k$, with equality if and only if $k\mid n$ and
$G\cong T_k(n)$. For every sufficiently large $n$, the unique maximizer
of $S_{\vartheta_k,k}$ is $T_k(n)$, and the graph with the second-largest value is
also unique. 
\end{theorem}

For $k=4$, the bound $S_{1,4}(G)\le3n/2$ follows from
Nikiforov~\cite[Theorem~5 and Corollary~6]{Nik2011}, who also records
$M_{1,4}(4t)=6t$ in~\cite[p.~2399]{Nik2011}. Theorem~\ref{main:order}
extends this determination to all sufficiently large orders, including
orders not divisible by $k$, and identifies the unique graph with the second-largest norm. It follows from Theorems~\ref{thm:critical}, \ref{thm:eventual}
and~\ref{thm:runnerup}. Corollary~\ref{cor:xi-four} gives the resulting
closed formula for $M_{1,4}(n)$ in each residue class modulo four.
Thus $k=4$ gives an answer to the classical Ky Fan problem for all
sufficiently large orders, whereas $k\ge5$ concerns the two-parameter norm
at $p=\vartheta_k$. We call $\vartheta_k$ a \emph{saturation exponent}
because the complete multipartite construction attains the order bound
at this exponent when $k\mid n$.

Nikiforov proved that $\xi_k(n)/n$ has a limit for every fixed $k$
\cite[Theorem~4]{Nik2011}. The same blowup argument gives the corresponding
limit for each fixed pair $(p,k)$; see Proposition~\ref{prop:asymptotic-limit}.
We also record a necessary arithmetic condition for finite-order equality
in Proposition~\ref{prop:arithmetic}. The distinction between finite-order
and asymptotic attainment is discussed in Section~\ref{sec:questions}.

The structural result used in the proof concerns the absolute deficit
$\Delta_k(G)=nC_k-S_{\vartheta_k,k}(G)$.
\begin{theorem}\label{main:structure}
For each fixed $k\ge4$, there are constants $\delta_k,K_k>0$ and $N_k$
such that, if $G$ has order $n\ge N_k$ and $\Delta_k(G)\le\delta_k$, there
is a complete $k$-partite graph $H=K_{n_1,\ldots,n_k}$ on the same vertex
set with
\[
 |E(G)\mathbin{\triangle}E(H)|+
 \sum_{i=1}^k(n_i-n/k)^2\le K_kn\Delta_k(G).
\]
For a sequence $|V(G_n)|=n\to\infty$, the condition $\Delta_k(G_n)=o(1)$
is equivalent to the existence of such graphs $H_n$ for which both the
edit count and the sum of squared part-size deviations are $o(n)$.
\end{theorem}

The structural aspect of Nikiforov's norm problems concerns not only the
extremal value, but also the extent to which a nearly maximal norm
determines the graph. The results above give a quantitative connection
between a small absolute norm deficit and multipartite structure.
This connection does not extend to relative near-optimality: connected
graphs may be asymptotically maximizing while remaining at quadratic
edge-edit distance from every complete $k$-partite graph. Thus uniqueness
of an extremal graph and stability of its near-extremal competitors are
distinct questions. This distinction complements the exact spectral
classifications studied in~\cite{Arevalo2021} and clarifies what
singular-value information can determine about graph structure.

We also consider how prescribing both the order and the size changes the
maximizing graphs. At $p=\vartheta_k$, a fixed deficit from
$|E(T_k(n))|$ that admits a complete $k$-partite realization confines all
maximizers to an explicit finite family of complete multipartite graphs;
unique extremizers are determined for specified deficits. A fixed
surplus gives a different extremal structure, in which additional edges
are placed inside a largest part of the balanced graph. These conclusions
hold over all graphs with the prescribed parameters, for fixed $k$ and
fixed edge deviation at sufficiently large orders. They describe how an
edge constraint affects the maximizing structure even when the leading
normalized norm is unchanged.

The two-parameter formulation retains both the exponent and the number
of singular values included in the norm. The upper estimates extend the
corresponding Schatten bounds to truncated spectra. On the lower-bound
side, a vertex-partition inequality complements the full-energy setting
of~\cite{Akbari2025}: it relates the norm of a graph to those of its
induced subgraphs under a prescribed total truncation index, with a best
universal remainder coefficient and complete equality conditions.
Together, the bounds and extremal descriptions connect Nikiforov's norm
questions with graph ordering, prescribed-size extremal problems and
near-extremal structure. They strengthen the link between matrix norm
inequalities and spectral extremal graph theory, and lead to further
questions about exponent ranges, finite-order attainment and truncated
minimum problems.

The paper is organized as follows. Section~\ref{sec:bounds} presents
general upper bounds, equality conditions and normalized extremal
limits. Section~\ref{sec:rigidity} concerns graph structure at small norm
deficits. Section~\ref{sec:local} studies stability and the two largest
norm values at prescribed order. Section~\ref{sec:prescribed} treats fixed
edge deficits and surpluses. Section~\ref{sec:partitions} gives the
vertex-partition lower bound and its equality cases. Further questions
are collected in Section~\ref{sec:questions}.

\section{General bounds and equality conditions}\label{sec:bounds}

Throughout, $n=|V(G)|$, $m=|E(G)|$, and $\rho=\lambda_1(G)=\sigma_1(G)$.
The equality $\rho=\sigma_1$ follows from nonnegativity of $A$.
We write $J_{a,b}$ for the all-ones matrix, $J_a=J_{a,a}$, and $I_a$ for the identity matrix.
For a real symmetric matrix $A$, write $|A|=(A^2)^{1/2}$. The unadorned matrix norm is the operator norm and $\|\cdot\|_F$ is the Frobenius norm.
All spectral lists count multiplicities. A secondary list obtained by removing a Perron value deletes exactly one occurrence of $\rho$, even when it is a repeated singular value.
For a matrix $B$, the notation $S_{p,k}(B)$ has the same meaning as
in~\eqref{eq:def}, with the singular values of $B$ in place of those of $G$.
Singular values beyond the smaller dimension are understood to be zero
when padding is convenient.
An assertion for all sufficiently large $n$ keeps the other named discrete parameters fixed; uniformity for a varying real exponent is stated explicitly.
We use the standard singular-value variational and interlacing facts collected in~\cite{HornJohnson2013}.

We shall use
\begin{equation}\label{eq:basics}
 \sum_{i=1}^n\sigma_i^2=\tr(A^2)=2m,
 \qquad \rho\ge\frac{2m}{n},\qquad \rho\le n-1.
\end{equation}
Equality in $\rho\ge2m/n$ holds if and only if $G$ is regular: the all-ones vector attains the largest Rayleigh quotient precisely in that case.

For the upper bounds in this section, assume $1\le p<2$ and $2\le k\le n$, and write
\[
 q=\frac p2,\qquad a=(k-1)^{1-q}.
\]
We first keep the number of edges in the estimate.

\begin{theorem}\label{thm:fixed-size}
Let $G$ be a nonempty graph with $n$ vertices and $m$ edges, and let $1\le p<2$ and $2\le k\le n$. Then
\begin{equation}\label{eq:rho-envelope}
 S_{p,k}(G)^p\le \rho^p+(k-1)^{1-p/2}(2m-\rho^2)^{p/2}.
\end{equation}
Put
\[
 \Gamma=\max\left\{\frac{2m}{n},\sqrt{\frac{2m}{k}}\right\}.
\]
Then
\begin{equation}\label{eq:fixed-envelope}
 S_{p,k}(G)^p\le \Gamma^p+(k-1)^{1-p/2}(2m-\Gamma^2)^{p/2}.
\end{equation}
Equality in \eqref{eq:fixed-envelope} holds if and only if
\[
 \rho=\Gamma,\qquad \sigma_2=\cdots=\sigma_k,
 \qquad \sigma_{k+1}=\cdots=\sigma_n=0.
\]
\end{theorem}
\begin{proof}
Strict concavity of $x^q$ gives
\[
 \sum_{i=2}^k\sigma_i^p
 \le (k-1)^{1-q}\left(\sum_{i=2}^k\sigma_i^2\right)^q
 \le a(2m-\rho^2)^q.
\]
Equality holds precisely when the selected $k-1$ squared singular values are equal and the omitted squared singular values vanish.
This proves \eqref{eq:rho-envelope} and its equality condition.

Let $f(x)=x^p+a(2m-x^2)^q$ for $0\le x\le\sqrt {2m}$.
On the interior,
\[
 f'(x)=p x\left[x^{p-2}-a(2m-x^2)^{q-1}\right].
\]
Because $q<1$, this derivative is positive for $x<\sqrt{2m/k}$, zero at $\sqrt{2m/k}$, and negative for $x>\sqrt{2m/k}$.
The largest possible value of $f$ on $[2m/n,\sqrt {2m}]$ is therefore $f(\Gamma)$.
Equation~\eqref{eq:basics} yields \eqref{eq:fixed-envelope}.
The maximizer of $f$ on this interval is unique, so equality also requires $\rho=\Gamma$.
\end{proof}

For an order-dependent bound, define
\[
F_{p,k}(t)=t^p+(k-1)^{1-p/2}[t(1-t)]^{p/2},\qquad 0\le t\le1.
\]

\begin{lemma}\label{lem:scalar}
The function $F_{p,k}$ has a unique maximizer $\tau_{p,k}\in(1/2,1)$.
Writing $z=(1-\tau_{p,k})/\tau_{p,k}$, this number is determined by
\begin{equation}\label{eq:z}
 (k-1)^{1-p/2}z^{p/2-1}(1-z)=2,\qquad 0<z<1.
\end{equation}
The maximum is nondegenerate: $F_{p,k}''(\tau_{p,k})<0$.
For fixed $k$, $\tau_{p,k}$ is strictly increasing in $p\in[1,2)$.
\end{lemma}
\begin{proof}
The derivative of $F_{p,k}$ is positive on $(0,1/2]$.
With $t=(1+z)^{-1}$, write
\[
 F_{p,k}(t)=\frac{1+a z^q}{(1+z)^p}.
\]
Differentiating with respect to $z$ gives
\[
 \frac{d}{dz}\frac{1+a z^q}{(1+z)^p}
 =q(1+z)^{-p-1}\left[a z^{q-1}(1-z)-2\right].
\]
The expression $g(z)=a z^{q-1}(1-z)$ is strictly decreasing on $(0,1)$, since
\[
 g'(z)=a z^{q-2}[(q-1)-qz]<0.
\]
It decreases from $+\infty$ to $0$, so \eqref{eq:z} has exactly one solution.
The derivative changes sign strictly, proving uniqueness. If $z_0$ is the
solution of~\eqref{eq:z}, then
\[
 \left.\frac{d^2}{dz^2}F_{p,k}\bigl((1+z)^{-1}\bigr)\right|_{z=z_0}
 =q(1+z_0)^{-p-1}g'(z_0)<0.
\]
Since $F'_{p,k}(\tau_{p,k})=0$ and $dt/dz=-(1+z)^{-2}\ne0$, the chain
rule gives $F''_{p,k}(\tau_{p,k})<0$.
Finally,
\[
 g(z)=\left(\frac{k-1}{z}\right)^{1-q}(1-z)
\]
strictly decreases as $q$ increases with $z$ fixed.
Its solution $g(z)=2$ therefore strictly decreases in $q$, and $\tau=(1+z)^{-1}$ strictly increases.
\end{proof}

Set
\[
C_{p,k}=F_{p,k}(\tau_{p,k})^{1/p}.
\]

\begin{theorem}\label{thm:order}
For $1\le p<2$, every graph $G$ of order $n\ge k\ge2$ satisfies
\begin{equation}\label{eq:order}
 S_{p,k}(G)\le nC_{p,k}.
\end{equation}
More precisely,
\begin{equation}\label{eq:finite-order}
 S_{p,k}(G)^p\le
 n^p F_{p,k}\!\left(\min\left\{\tau_{p,k},1-\frac1n\right\}\right).
\end{equation}
For a nonempty graph, equality in \eqref{eq:order} holds if and only if $G$ is regular and
\[
 \rho=n\tau_{p,k},\qquad
 \sigma_2=\cdots=\sigma_k
 =n\sqrt{\frac{\tau_{p,k}(1-\tau_{p,k})}{k-1}},
 \qquad \sigma_i=0\quad(i>k).
\]
\end{theorem}
\begin{proof}
The empty graph satisfies both inequalities strictly.
For a nonempty graph, \eqref{eq:rho-envelope} and $2m\le n\rho$ imply
\[
 S_{p,k}(G)^p
 \le\rho^p+a(2m-\rho^2)^q
 \le\rho^p+a(n\rho-\rho^2)^q
 =n^pF_{p,k}(\rho/n).
\]
Now use Lemma~\ref{lem:scalar} and $\rho/n\le1-1/n$.
For equality in \eqref{eq:order}, strict monotonicity in $2m$ requires $2m=n\rho$, which is equivalent to regularity.
The remaining conditions are precisely the concavity equality condition and the unique scalar maximizer.
Conversely, the stated spectrum and regularity make every inequality an equality.
\end{proof}

\begin{proposition}\label{prop:comparison}
For $k\ge2$,
\[
 C_{1,k}=\frac{1+\sqrt{k}}2.
\]
For $1<p<2$,
\begin{equation}\label{eq:strict-improvement}
 C_{p,k}<\frac{1+k^{1/p-1/2}}2.
\end{equation}
\end{proposition}
\begin{proof}
For $p=1$, maximizing $t+\sqrt{(k-1)t(1-t)}$ gives
$\tau=(1+k^{-1/2})/2$ and the stated value.
For $p>1$, let $t=\tau_{p,k}$ and $s=\sqrt{t(1-t)/(k-1)}$.
For the vector $x=(t,s,\ldots,s)\in\R^k$, write $x=e_1/2+v$.
Then
\[
 \|v\|_2^2=(t-1/2)^2+(k-1)s^2=\frac14.
\]
Since $s>0$, the two summands are not collinear.
Strictness of the triangle inequality in $\ell_p^k$ and the power-mean inequality give
\[
 C_{p,k}=\|x\|_p<\frac12+\|v\|_p
 \le\frac12+\frac12 k^{1/p-1/2}.
\]
This completes the proof.
\end{proof}

\begin{remark}
The right-hand side of \eqref{eq:strict-improvement} is the direct bound obtained from $2A=J_n+(2A-J_n)$, the triangle inequality, and the Frobenius norm of $2A-J_n$.
At $p=1$, Theorem~\ref{thm:order} recovers Nikiforov's classical bound~\cite{Nik2011}.
For $1<p<2$, the improvement is strict, but this comparison alone does not assert that $C_{p,k}$ is attainable for every parameter pair.
\end{remark}

We next identify an exponent at which the order bound is attained.
When $k\mid n$, all parts of $T_k(n)$ have size $n/k$.
The spectral classification used in the equality case is classical
\cite[Lemma~2.3 and Theorem~2.5]{Arevalo2021}; a direct argument is
included below.

\begin{theorem}\label{thm:critical}
Let $k\ge4$ and let $\vartheta_k$ be given by \eqref{eq:critical-exponent}.
Every graph $G$ of order $n\ge k$ satisfies
\begin{equation}\label{eq:main}
 S_{\vartheta_k,k}(G)\le \frac{n}{k}\binom{k}{2}^{1/\vartheta_k}.
\end{equation}
Equality holds if and only if $k\mid n$ and $G\cong T_k(n)$.
\end{theorem}
\begin{proof}
Set $p=\vartheta_k$.
The identity
\[
 (k-1)^p=\frac{(k-1)(k-2)}2
\]
shows that $z=1/(k-1)$ satisfies \eqref{eq:z}.
Consequently,
\[
 \tau_{p,k}=\frac{k-1}{k},\qquad
 C_{p,k}^p=\frac{(k-1)^p+(k-1)}{k^p}
          =\frac{\binom{k}{2}}{k^p}.
\]
Inequality~\eqref{eq:main} follows from Theorem~\ref{thm:order}.

Suppose equality holds.
The graph is regular of degree $d=(k-1)n/k$, has rank $k$, and has secondary nonzero singular values $s=n/k$.
Since $d>s$, the singular value $d$ is simple.
A disconnected $d$-regular graph would have the eigenvalue $d$ at least twice, so $G$ is connected.
All $k-1$ other nonzero eigenvalues are $s$ or $-s$.
The equality $\tr A=0$ and $d=(k-1)s$ force all of them to be $-s$.
Therefore $A$ is negative semidefinite on $\one^\perp$.
The complete multipartite characterization
\cite[Lemma~2.3]{Arevalo2021} now applies. We include the argument.

Let $u,v$ be nonadjacent vertices and put $x=e_u-e_v$.
Then $x\perp\one$ and $x^TAx=0$.
The spectral decomposition of the negative-semidefinite restriction implies $Ax=0$.
Thus $u$ and $v$ have identical open neighborhoods.
The relation ``equal or nonadjacent'' is now an equivalence relation.
Indeed, if distinct $u,v,w$ satisfy $uv,vw\notin E(G)$, then $N(u)=N(v)=N(w)$; an edge $uw$ would imply $w\in N(w)$, which is impossible.
Each equivalence class is independent and any two distinct classes are completely joined.
Hence $G$ is complete multipartite.

Regularity forces all part sizes to be equal.
A complete $r$-partite graph with nonempty parts has rank $r$: its adjacency matrix on the normalized part indicators is congruent to $J_r-I_r$, which is nonsingular for $r\ge2$, and it vanishes on the within-part zero-sum subspaces.
Since $\rank A=k$, there are exactly $k$ parts.
Thus $k\mid n$ and $G\cong T_k(n)$.

Conversely, for $n=kt$, the spectrum of $T_k(n)$ is
\[
 (k-1)t,\quad (-t)^{[k-1]},\quad 0^{[k(t-1)]}.
\]
Substitution in \eqref{eq:def} gives equality.
\end{proof}

\subsection{Exact equality and normalized extrema}

\begin{proposition}\label{prop:arithmetic}
Suppose that $S_{p,k}(G)=nC_{p,k}$ for a nonempty graph, where
$1\le p<2$ and $k\ge2$. Then an integer $h$ satisfies
\[
 1<h\le k-1,\qquad h\equiv k-1\pmod2,\qquad
 2h^p=h^2-(k-1),
 \qquad \tau_{p,k}=\frac{h^2}{h^2+k-1}.
\]
For fixed $k$, exact equality is therefore possible at only finitely many
exponents in $[1,2)$. If $p=1$, necessarily $k=(h-1)^2$.
\end{proposition}
\begin{proof}
By Theorem~\ref{thm:order}, the graph is regular of degree
$d=n\tau_{p,k}$, has rank $k$, and its remaining nonzero singular values
are equal to a number $s>0$. The moment identity gives
$d^2/s^2=(k-1)\tau_{p,k}/(1-\tau_{p,k})>k-1$, so $d>s$.
If $u$ and $v$ count the eigenvalues $s$ and $-s$, respectively, then
$u+v=k-1$ and zero trace gives $d=(v-u)s$. Set $h=v-u$.
Thus $h>1$, with the stated upper bound and parity. From
$nd=d^2+(k-1)s^2$ we obtain the formula for $\tau_{p,k}$.
Substituting $z=(k-1)/h^2$ into~\eqref{eq:z} gives
$h^{-p}(h^2-k+1)=2$. Each of the finitely many integers $h>1$
determines at most one exponent, and $p=1$ gives $k=(h-1)^2$.
\end{proof}

The next statement extends \cite[Theorem~4]{Nik2011} by the same
blowup argument.
\begin{proposition}\label{prop:asymptotic-limit}
For every fixed integer $k\ge1$ and real $1\le p<\infty$, the limit
\[
 \mathfrak m_{p,k}:=\lim_{n\to\infty}\frac{M_{p,k}(n)}n
 =\sup_{h\ge k}\frac{M_{p,k}(h)}h
\]
exists. If $1\le p\le q<\infty$, then
\[
 \mathfrak m_{q,k}\le\mathfrak m_{p,k}
 \le k^{1/p-1/q}\mathfrak m_{q,k}.
\]
In particular, $p\mapsto\mathfrak m_{p,k}$ is continuous,
$\mathfrak m_{p,k}=1$ for $p\ge2$, and
$\mathfrak m_{\vartheta_k,k}=C_k$ for $k\ge4$.
\end{proposition}
\begin{proof}
The supremum is finite, since $S_{p,k}(G)\le k^{1/p}(n-1)$.
Fix $h\ge k$ and an $h$-vertex graph $H$ attaining $M_{p,k}(h)$.
Replacing each vertex by an independent set of size $t$ gives an adjacency
matrix $A(H)\otimes J_t$, so all nonzero singular values are multiplied
by $t$. For $n=ht+r$, $0\le r<h$, add $r$ isolated vertices. Therefore
$M_{p,k}(n)\ge tM_{p,k}(h)$, and
\[
 \liminf_{n\to\infty}\frac{M_{p,k}(n)}n
 \ge\frac{M_{p,k}(h)}h.
\]
Taking the supremum over $h$ and comparing with the limsup proves
existence. The displayed comparison follows from the corresponding
$\ell_p$ inequalities for the first $k$ singular values; it also proves
continuity. For $p\ge2$,
$S_{p,k}(G)\le S_{2,n}(G)\le\sqrt{n(n-1)}$, whereas
$S_{p,k}(K_n)\ge n-1$. Finally, Theorem~\ref{thm:critical} bounds the
last limit by $C_k$ and attains it along $n=kt$.
\end{proof}

\section{Norm deficits and multipartite structure}\label{sec:rigidity}

Retain $1\le p<2$, $k\ge2$, and $q=p/2$.
For a graph of order $n\ge k$, set
\[
 R_k=\sum_{i>k}\sigma_i^2,\qquad
 b_k=\frac1{k-1}\sum_{i=2}^k\sigma_i^2,
 \qquad \mathcal D=n^pC_{p,k}^p-S_{p,k}(G)^p.
\]
Put
\[
 \alpha_{p,k}=q[4(k-1)]^{1-q},\qquad
 \beta_p=\frac{q(1-q)}2\,4^{2-q}.
\]

\begin{theorem}\label{thm:defect}
With the above notation,
\begin{align}\label{eq:defect}
 \mathcal D\ \ge\;& n^p\left[F_{p,k}(\tau_{p,k})-F_{p,k}(\rho/n)\right]\notag\\
 &+\alpha_{p,k}n^{p-2}(n\rho-2m+R_k)\notag\\
 &+\beta_p n^{p-4}\sum_{i=2}^k(\sigma_i^2-b_k)^2.
\end{align}
Every term on the right is nonnegative.
\end{theorem}
\begin{proof}
The empty graph follows directly, so assume $m>0$.
Let
\[
 U=n\rho-\rho^2,\qquad
 T=\sum_{i=2}^k\sigma_i^2=2m-\rho^2-R_k.
\]
Then $0\le T\le U\le n^2/4$, and
\[
 U-T=n\rho-2m+R_k.
\]
Since $x^{q-1}$ is decreasing,
\[
 U^q-T^q=\int_T^U qx^{q-1}\,dx
 \ge q(n^2/4)^{q-1}(U-T).
\]
Multiplying by $a=(k-1)^{1-q}$ gives the second term in \eqref{eq:defect}.

For $0<x\le n^2/4$,
\[
 -\frac{d^2}{dx^2}x^q=q(1-q)x^{q-2}
 \ge q(1-q)4^{2-q}n^{p-4}.
\]
Thus $x^q+\beta_p n^{p-4}x^2$ is concave on $[0,n^2/4]$, by continuous extension at zero. Each selected square satisfies $\sigma_i^2\le T\le n^2/4$.
Apply Jensen's inequality to $x_i=\sigma_i^2$, $2\le i\le k$, to obtain
\[
 aT^q-\sum_{i=2}^k\sigma_i^p
 \ge\beta_p n^{p-4}\sum_{i=2}^k(x_i-b_k)^2.
\]
Finally use the exact decomposition
\begin{align*}
 \mathcal D={}&n^p\left[F_{p,k}(\tau_{p,k})-F_{p,k}(\rho/n)\right]\\
 &+a(U^q-T^q)+\left(aT^q-\sum_{i=2}^k\sigma_i^p\right).
\end{align*}
This proves the claim.
\end{proof}

\subsection{Binary matrices close to a scaled projection}

\begin{lemma}\label{lem:rounding}
For every integer $k\ge1$ there is $\varepsilon_k>0$ with the following property.
Let $n\ge k$, let $R$ be a symmetric $(0,1)$ matrix with diagonal entries one, and let $P$ be an orthogonal projection of rank $k$.
Put
\[
 a=n/k,\qquad D=\Frob{R-aP},\qquad r_0=\max_v|(R\one)_v-a|.
\]
If $D^2\le\varepsilon_k n$ and $r_0\le\varepsilon_k n$, there is a partition into $k$ nonempty classes, with associated all-ones block matrix $B$, such that
\begin{equation}\label{eq:scalable-rounding}
 \Frob{R-B}\le2D,\qquad
 \sum_{i=1}^k(n_i-n/k)^2\le9D^2.
\end{equation}
In particular, the conclusion holds for all sufficiently large $n$ if $D$ is bounded and $r_0=O(\sqrt n)$.
\end{lemma}
\begin{proof}
The exact classification of nonnegative symmetric idempotents gives positive rank-one blocks~\cite[Chapter~3, Corollary~3.4]{BermanPlemmons1994}.
We prove the approximate binary statement, including its dependence on the error.
Set
\[
 c=(2k)^{-(2k-1)/2},\qquad
 \eta=\min\{(64k^2)^{-1},c^2/16\},
\]
and choose
\[
 0<\varepsilon_k\le
 \min\left\{\frac{1}{2k(\eta^{-2}+\eta^{-1})},
             \frac{1}{4k(5k+9)}\right\}.
\]
Write $M=aP=ZZ^T$, where $Z^TZ=aI_k$, and denote the rows of $Z$ by $z_v$.
Call a row good if $\|(R-M)_{v,*}\|_2\le\eta$ and let $b_0$ be the number of other rows. Then $b_0\le D^2/\eta^2$.
The diagonal condition and Cauchy's inequality give the crucial aggregate estimate
\[
 \sum_{v\text{ not good}}\|z_v\|_2^2
 =\sum_{v\text{ not good}}M_{vv}
 \le b_0+\sqrt{b_0}D
 \le D^2(\eta^{-2}+\eta^{-1})\le\frac{n}{2k}.
\]
Consequently
\begin{equation}\label{eq:good-covariance}
 \sum_{v\text{ good}}z_v^Tz_v\succeq\frac{n}{2k}I_k.
\end{equation}
Every good row has squared length at most $1+\eta\le2$.
Suppose that the previously selected rows span a proper subspace $L$ of
$\R^k$, and let $\Pi_{L^\perp}$ be the orthogonal projection onto its
orthogonal complement. Taking the trace after compression
in~\eqref{eq:good-covariance} gives
\[
 \sum_{v\text{ good}}\|z_v\Pi_{L^\perp}\|_2^2
 =\tr\!\left(\Pi_{L^\perp}\sum_{v\text{ good}}z_v^Tz_v\right)
 \ge\frac{n}{2k}\dim L^\perp.
\]
There are at most $n$ good rows. Hence one of them has squared distance
at least $1/(2k)$ from $L$. Repeating this choice yields $k$ good rows.
Their row matrix $Y$ satisfies $|\det Y|\ge(2k)^{-k/2}$ and $\|Y\|\le\sqrt{2k}$, hence $\sigma_{\min}(Y)\ge c$.
Write $i_1,\ldots,i_k$ for their indices.

If $R_{i_j i_\ell}=1$ for $j\ne\ell$, the entrywise errors in the selected rows imply
\[
 \|z_{i_j}-z_{i_\ell}\|_2^2
 =M_{i_ji_j}+M_{i_\ell i_\ell}-2M_{i_ji_\ell}\le4\eta.
\]
This contradicts $\|(e_j-e_\ell)^TY\|_2^2\ge2c^2$.
Thus $R[I,I]=I_k$ and $\|YY^T-I_k\|\le k\eta$.
For every vertex put
\[
 h_v=(R_{v i_1},\ldots,R_{v i_k})^T\in\{0,1\}^k.
\]
Symmetry and goodness of the selected rows give $Yz_v^T=h_v+e_v$, where $\|e_v\|_2\le\sqrt{k}\eta$ uniformly in $v$.
Since $\|(YY^T)^{-1}\|\le2$ and $\|(YY^T)^{-1}-I_k\|\le2k\eta$, we obtain
\begin{equation}\label{eq:binary-pattern-gram}
 |M_{vw}-h_v^Th_w|\le8k^2\eta\le\tfrac14\qquad(v,w\in[n]).
\end{equation}
Indeed, expand $(h_v+e_v)^T(YY^T)^{-1}(h_w+e_w)$ and use $\|h_v\|_2\le\sqrt k$.

Call a vertex exceptional if $\|h_v\|_1\ne1$.
At such a vertex, \eqref{eq:binary-pattern-gram} and $R_{vv}=1$ imply $|R_{vv}-M_{vv}|\ge3/4$, so there are at most $4D^2$ exceptional vertices.
Let $N_i$ count the nonexceptional vertices with pattern $e_i$.
The selected row $i_i$ then gives
\[
 |N_i-a|\le r_0+4D^2.
\]
For a fixed $v$, every disagreement $R_{vw}\ne h_{v,i}$ at a nonexceptional vertex of pattern $e_i$ contributes at least $1/4$ to the squared error of row $v$.
There are at most $4D^2$ such disagreements, and at most $4D^2$ exceptional columns.
It follows that
\[
 \left|(R\one)_v-\sum_i h_{v,i}N_i\right|\le8D^2.
\]
Writing $\ell_v=\|h_v\|_1\le k$ yields
\[
 |(\ell_v-1)a|
 \le(k+1)r_0+(4k+8)D^2
 \le(5k+9)\varepsilon_k n<a.
\]
Since $\ell_v$ is an integer, $\ell_v=1$ for every $v$.
Partition by these unit patterns. All $k$ classes are nonempty, since the anchor $i_i$ has pattern $e_i$.
The associated matrix is $B_{vw}=h_v^Th_w$.
Every nonzero entry of $R-B$ contributes at least $1/4$ to $\Frob{R-M}^2$, by \eqref{eq:binary-pattern-gram}; hence $\Frob{R-B}\le2D$ and $\Frob{B-M}\le3D$.
The nonzero eigenvalues of $B$ are the class sizes and those of $M$ are $a$ repeated $k$ times.
The Hermitian Hoffman--Wielandt inequality~\cite[Corollary~6.3.8]{HornJohnson2013} now gives the second assertion in \eqref{eq:scalable-rounding}.
\end{proof}

\subsection{Recovering the signs and the graph partition}

\begin{proposition}\label{prop:spectral-criterion}
Fix $k\ge4$, and for a graph of order $n\ge k$ put
\begin{equation}\label{eq:spectral-defect}
\begin{split}
 \mathfrak r_k(G)={}&n\rho-2m+\sum_{i>k}\sigma_i^2
 +\left(\rho-\frac{k-1}{k}n\right)^2
 +\sum_{i=2}^k\left(\sigma_i-\frac nk\right)^2.
\end{split}
\end{equation}
There are constants $\varepsilon,K>0$ and $N$, depending only on $k$, such that
$\mathfrak r_k(G)\le\varepsilon n$ and $n\ge N$ imply the existence of a complete
$k$-partite graph $H=K_{n_1,\ldots,n_k}$ on the same vertex set satisfying
\[
 |E(G)\mathbin{\triangle}E(H)|+\sum_{i=1}^k(n_i-n/k)^2
 \le K\mathfrak r_k(G).
\]
\end{proposition}
\begin{proof}
All constants in this proof depend only on $k$. Set $r=\mathfrak r_k(G)$,
$a=n/k$ and $b=(k-1)/k$. Each term in~\eqref{eq:spectral-defect} is nonnegative.
In particular,
\[
 |\rho-bn|\le\sqrt r,\qquad |\sigma_i-a|\le\sqrt r\quad(2\le i\le k),
 \qquad \sum_{i>k}\sigma_i^2\le r.
\]
Choose eigenvalues corresponding to the first $k$ singular values, including
one Perron eigenvalue. If $j\ge1$ of the selected secondary eigenvalues were
positive, their signed sum, including the Perron value, would be at least
$2jn/k-k\sqrt r$. The absolute sum of the omitted eigenvalues is at most
$\sqrt{nr}$. Choose $\varepsilon$ so that $(k+1)\sqrt\varepsilon<2/k$.
Since $r\le\varepsilon n$ and $n\ge1$, these estimates are incompatible with
$\tr A=0$. Thus all selected secondary eigenvalues are negative.
After increasing $N$, the selected magnitudes are of order $n$ and the omitted
ones are at most $\sqrt r=O(\sqrt n)$; hence the selected spectral subspace is
well defined. In particular, $\rho$ is simple and
$\rho-\lambda_2(A)\ge c_kn$ for some $c_k>0$.

Let $u$ be a unit nonnegative Perron vector, $w=n^{-1/2}\one$, and let $P$
project onto the selected spectral subspace. The Rayleigh quotient gives
\[
 1-|u^Tw|^2\le\frac{\rho-w^TAw}{\rho-\lambda_2(A)}
 =\frac{n\rho-2m}{n(\rho-\lambda_2(A))}\le C r/n^2.
\]
Therefore $\|nuu^T-J_n\|_F^2\le C r$. Moreover, by the signs just proved,
\[
 \|A-(nuu^T-aP)\|_F^2
 = (\rho-bn)^2+\sum_{i=2}^k(\sigma_i-a)^2+\sum_{i>k}\sigma_i^2
 \le r.
\]
It follows that $D^2:=\|(J_n-A)-aP\|_F^2\le C r$.

Put $\bar d=2m/n$. Since
$\sum_v d(v)^2=\|A\one\|_2^2\le n\rho^2$ and
$\rho+\bar d\le2n$, the degree variance satisfies
\[
 \begin{aligned}
 \sum_v(d(v)-\bar d)^2
 &\le n(\rho^2-\bar d^{\,2})
   =(n\rho-2m)(\rho+\bar d)\\
 &\le2n(n\rho-2m)\le2nr.
 \end{aligned}
\]
Moreover, $|\bar d-bn|\le\sqrt r+r/n$.
Consequently the row sums of $R=J_n-A$ differ from $a$ by at most
$\sqrt{2nr}+\sqrt r+r/n\le C\sqrt{nr}$, after taking $\varepsilon\le1$.
This is a symmetric binary matrix with diagonal one. Take $\varepsilon$
small enough that $C\varepsilon\le\varepsilon_k$ and
$C\sqrt\varepsilon\le\varepsilon_k$, where $\varepsilon_k$ is the
constant in Lemma~\ref{lem:rounding}. Then $r\le\varepsilon n$ implies
both hypotheses of that lemma.
That lemma gives an all-ones block matrix $B$ with
\[
 \tfrac12\|R-B\|_F^2\le2D^2,\qquad
 \sum_i(n_i-a)^2\le9D^2.
\]
The graph with adjacency matrix $J_n-B$ is the required $H$.
\end{proof}

\begin{theorem}\label{thm:multiscale}
Fix $k\ge4$. There are $\delta_0,K>0$ and $N$, depending only on $k$, such
that every graph $G$ of order $n\ge N$ with
$\Delta_k(G)=nC_k-S_{\vartheta_k,k}(G)\le\delta_0$ admits a complete $k$-partite
graph $H=K_{n_1,\ldots,n_k}$ on the same vertex set satisfying
\begin{equation}\label{eq:multiscale}
 |E(G)\mathbin{\triangle}E(H)|+
 \sum_{i=1}^k(n_i-n/k)^2\le Kn\Delta_k(G).
\end{equation}
\end{theorem}
\begin{proof}
Put $p=\vartheta_k$, $b=(k-1)/k$ and $\Delta=\Delta_k(G)$. All constants depend
only on $k$. By Lemma~\ref{lem:scalar}, uniqueness and nondegeneracy of the
scalar maximum give, on the whole interval $[0,1]$,
\[
 F_{p,k}(b)-F_{p,k}(t)\ge c_0(t-b)^2
\]
for some $c_0>0$. The mean value theorem gives
$\mathcal D\le Cn^{p-1}\Delta$. Theorem~\ref{thm:defect} therefore yields
\[
 R_k+n\rho-2m\le Cn\Delta,\qquad
 |\rho-bn|\le C\sqrt{n\Delta},\qquad
 \sum_{i=2}^k(\sigma_i^2-b_k)^2\le Cn^3\Delta.
\]
Since
\[
 b_k=\frac{n\rho-\rho^2-(n\rho-2m+R_k)}{k-1},
\]
we have $|b_k-n^2/k^2|\le Cn^{3/2}\sqrt\Delta$ once $\Delta\le\delta_0$
and $n$ is sufficiently large. Consequently
$|\sigma_i^2-n^2/k^2|\le Cn^{3/2}\sqrt\Delta$ for $2\le i\le k$.
Dividing by $\sigma_i+n/k\ge n/k$ gives
$|\sigma_i-n/k|\le C\sqrt{n\Delta}$. Thus
\[
 \mathfrak r_k(G)\le Cn\Delta.
\]
Choose $\delta_0$ so that $C\delta_0$ is below the constant in
Proposition~\ref{prop:spectral-criterion}. That proposition proves
\eqref{eq:multiscale}.
\end{proof}

\begin{corollary}\label{thm:bounded-edit}
Fix $k\ge4$ and $K_0\ge0$.
For all sufficiently large $n$, every graph with
$0\le nC_k-S_{\vartheta_k,k}(G)\le K_0/n$ differs in at most $B_0$ edges from a complete $k$-partite graph $K_{n_1,\ldots,n_k}$ with $|n_i-n/k|\le B_0$.
The constants $B_0$ and the order threshold depend only on $k,K_0$.
\end{corollary}
\begin{proof}
Apply Theorem~\ref{thm:multiscale} at $p=\vartheta_k$.
The sum of the edit count and the squared part deviations is at most $KK_0$.
Take $B_0\ge\max\{KK_0,\sqrt{KK_0}\}$ and then enlarge the order threshold so that $K_0/n\le\delta_0$.
\end{proof}

In the reciprocal case the preceding proof also gives
\[
 \rho=\frac{k-1}{k}n+O(1),\qquad
 \sigma_i=\frac nk+O(1)\ (2\le i\le k),\qquad R_k+n\rho-2m=O(1).
\]
The converse edit statement will follow from the local expansion, proved independently in the next section.

\section{Local losses, stability and extremal graphs of prescribed order}\label{sec:local}
Fix $k\ge4$ and put $\kappa=k-1$.
For $1\le p<2$, define
\begin{align*}
 D_p&=(\kappa^p+\kappa)^{1/p},&c_T(p)&=D_p/k,\\
 g_k(p)&=\frac{D_p^{1-p}}{2k}
 \{2\kappa(\kappa^{p-1}+1)-k(p-1)(k-2)\},\\
 z_+(p)&=2D_p^{1-p}(\kappa-\kappa^{p-1}),&
 z_-(p)&=2D_p^{1-p}(\kappa^{p-1}+1).
\end{align*}
At $p=\vartheta_k$, set
\[
 \gamma_k=g_k(\vartheta_k)=\tfrac12(kC_k)^{1-\vartheta_k}[1+(2-\vartheta_k)(k-2)],
 \qquad \zeta_k=z_+(\vartheta_k)=z_-(\vartheta_k)=k(kC_k)^{1-\vartheta_k}.
\]
Both are positive, and
\[
 \zeta_k-2\gamma_k=(kC_k)^{1-\vartheta_k}[1+(\vartheta_k-1)(k-2)]>0.
\]
Choose a compact neighborhood $J\subset[1,2)$ of $\vartheta_k$, one-sided at $\vartheta_4=1$, on which $g_k,z_+,z_-$ have a common positive lower bound.
Every estimate indexed by $k,J$ is uniform over $p\in J$ and the explicitly indicated graphs.

\subsection{Part-size perturbations}
\begin{lemma}\label{lem:secondorder}\label{thm:cubic}
Let $H=K_{n_1,\ldots,n_k}$, $\delta_i=n_i-n/k$, and $V=\sum_i\delta_i^2$.
Uniformly for $p\in J$, $V\le n$, and sufficiently large $n$,
\begin{equation}\label{eq:general-part-loss}
 S_{p,k}(H)=nc_T(p)-\frac{g_k(p)V}{n}
             +O_{k,J}(V^{3/2}/n^2).
\end{equation}
If $|\delta_i|\le B_1$ and $p=\vartheta_k$, then uniformly
\begin{equation}\label{eq:cubic-expansion}
 S_{\vartheta_k,k}(H)=nC_k-\frac{\gamma_k}{n}\sum_i\delta_i^2
 +\frac{\xi_k}{n^2}\sum_i\delta_i^3+O_{k,B_1}(n^{-3}),
\end{equation}
where
\begin{equation}\label{eq:xi}
 \xi_k=(kC_k)^{1-\vartheta_k}
 \left[(k-1)\left(\frac32-\vartheta_k\right)
 +\frac{k(k-3)}6(\vartheta_k-1)(\vartheta_k-2)\right].
\end{equation}
Moreover, $\xi_4,\xi_5>0$, while $\xi_k<0$ for every $k\ge6$.
\end{lemma}
\begin{proof}
For $x_i=n_i/n$, the nonzero eigenvalues of $A(H)/n$ are those of
$Q(x)=\sqrt{x}\sqrt{x}^{\,T}-\diag x$; this follows by restricting to the normalized part indicators.
At $x=k^{-1}\one$ the eigenvalues are $(k-1)/k$ and $-1/k$, the latter repeated $k-1$ times.
Consequently $\Phi_p(x)=S_{p,k}(Q(x))$ is smooth jointly in $p,x$ near that point, including at $p=1$.
For example, separated contours around the positive and negative eigenvalues define analytic branches of $z^p$ and $(-z)^p$ and hence their trace functions.
There is no differentiability assertion at a zero eigenvalue.
Permutation symmetry makes the linear term vanish on $\one^\perp$, the quadratic term a multiple of $\sum_i h_i^2$, and the cubic term a multiple of $\sum_i h_i^3$.
For the last assertion, express a symmetric cubic using the monomial types $h_i^3$, $h_i^2h_j$, and $h_ih_jh_l$, and use $\sum_i h_i=0$.

It suffices to evaluate $x_1=1/k+(k-1)u$ and $x_2=\cdots=x_k=y=1/k-u$.
There are $k-2$ eigenvalues $-y$; the remaining positive and negative eigenvalues are $\lambda_+(u)$ and $-\nu(u)$.
The discriminant of their two-dimensional quotient is
\[
 1+2(k-2)u-k(3k-4)u^2.
\]
Writing $a=1/k$, $b=(k-1)/k$, and $\kappa=k-1$, its square-root expansion gives
\begin{align*}
 \lambda_+(u)&=b-\kappa^2u^2+(k-2)\kappa^2u^3+O_k(u^4),\\
 \nu(u)&=a+(k-2)u-\kappa^2u^2+(k-2)\kappa^2u^3+O_k(u^4).
\end{align*}
The quadratic and cubic coefficients of
$\lambda_+(u)^p+\nu(u)^p+(k-2)y^p$ are respectively
\begin{align*}
 A_2={}&-p\kappa^2(b^{p-1}+a^{p-1})
       +\frac{p(p-1)}2a^{p-2}(k-2)\kappa,\\
 A_3={}&p(k-2)\kappa^2\{b^{p-1}+a^{p-1}-(p-1)a^{p-2}\}\\
 &+\frac{p(p-1)(p-2)}6a^{p-3}\{(k-2)^3-(k-2)\}.
\end{align*}
There is no linear term.
Taking the $p$th root multiplies these coefficients by $c_T(p)^{1-p}/p$.
The chosen direction has squared norm $k\kappa$ and cubic sum $k\kappa(k-2)$.
Therefore the quadratic coefficient of $\Phi_p$ on $\one^\perp$ is $-g_k(p)$.
At $p=\vartheta_k$, the identity $\kappa^{p-1}=(k-2)/2$ gives the cubic coefficient $\xi_k$ in \eqref{eq:xi}.

Taylor's theorem on a fixed compact neighborhood supplies a uniform cubic remainder for the quadratic expansion.
Substituting $x_i=1/k+\delta_i/n$ and multiplying by $n$ gives \eqref{eq:general-part-loss} when $V\le n$.
A uniform fourth-order remainder gives \eqref{eq:cubic-expansion} for bounded deviations.
For the signs, drop the positive factor in \eqref{eq:xi}.
At $k=4$ the bracket is $3/2$.
For $k=5$, $1<\vartheta_5<4/3$ follows from $6^3<4^4$; the bracket
$4(3/2-p)+(5/3)(p-1)(p-2)$ decreases on $[1,4/3]$ and equals $8/27$ at $4/3$.
For $k=6$, $7/5<\vartheta_6<3/2$ follows from $5^7<10^5$ and $10^2<5^3$; the bracket
$5(3/2-p)+3(p-1)(p-2)$ decreases on $[7/5,3/2]$ and equals $-11/50$ at $7/5$.
Finally, $k\ge7$ implies $k-2\ge2\sqrt{k-1}$, so $3/2\le \vartheta_k<2$ and the bracket is negative.
\end{proof}
For $n=kt+r$, write
\[
 v_r=\frac{r(k-r)}k,\qquad w_r=\frac{r(k-r)(k-2r)}{k^2}.
\]
The balanced part deviations have square sum $v_r$ and cubic sum $w_r$; hence
\begin{equation}\label{eq:balanced-expansion}
 S_{\vartheta_k,k}(T_k(n))=nC_k-\frac{\gamma_kv_r}{n}
                  +\frac{\xi_kw_r}{n^2}+O_k(n^{-3}).
\end{equation}

\subsection{A localized estimate for a separated spectral cluster}
\begin{proposition}\label{prop:localized}\label{lem:finite-edges}
Fix $k\ge4$ and the compact interval $J$ above.
Let $H=K_{n_1,\ldots,n_k}$ have $V\le n$, and let a symmetric edit matrix $E$ add or delete $e\le n$ graph edges.
Let $P_0$ project onto the normalized part-indicator subspace. For the
edit graph, write $d_E(v)=|\{w:E_{vw}\ne0\}|$ and
$d_E=\max_v d_E(v)$.
For all sufficiently large $n$, uniformly in $p\in J$ and $0\le t\le1$, the matrix $A_t=A(H)+tE$ has exactly $k$ eigenvalues of magnitude at least $c_kn$, while the others have magnitude at most $\sqrt{2e}$.
If $P_t$ is the projection onto that group, then
\begin{equation}\label{eq:localized-E}
 \|EP_0\|_F^2\le\frac{4ked_E}{n},\qquad
 \|EP_t\|_F\le2\|EP_0\|_F,
\end{equation}
and
\begin{equation}\label{eq:localized-second}
 \left|\frac{d^2}{dt^2}S_{p,k}(A_t)\right|
 \le\frac{C_{k,J}}n\|EP_t\|_F^2
 \le C_{k,J}\frac{ed_E}{n^2}.
\end{equation}
The constants and the order threshold do not depend on the graph, the positions of the edited edges, or $p\in J$.
\end{proposition}
\begin{proof}
Fix $k$ and $J$. All constants below are independent of $H$, of the
positions of the edited edges, of $t\in[0,1]$ and of $p\in J$.
Write $A_0=A(H)$ and let $U$ have the normalized part indicators as its
columns, so $P_0=UU^T$. The assumption $V\le n$ gives
$n_i\ge n/k-\sqrt n\ge n/(2k)$ once $n$ is sufficiently large. Hence
\[
\begin{aligned}
 \|E\|&\le\|E\|_F=\sqrt{2e},\\
 \|EP_0\|_F^2&=\|EU\|_F^2
 =\sum_{v,i}\frac{\bigl(\sum_{w\in V_i}E_{vw}\bigr)^2}{n_i}\\
 &\le\frac{2k}{n}\sum_v d_E(v)^2
 \le\frac{4ked_E}{n}.
\end{aligned}
\]

We first choose an order threshold for spectral separation. Put
$x_i=n_i/n$. On the range of $P_0$, the matrix $A_0/n$ is represented by
$Q(x)=\sqrt{x}\sqrt{x}^{\,T}-\diag x$, and it is zero on the orthogonal
complement. Since $\|x-k^{-1}\one\|_2\le n^{-1/2}$,
\[
 \left\|Q(x)-\frac{J_k-I_k}{k}\right\|\le C_k n^{-1/2}.
\]
Also $\|tE\|/n\le\sqrt{2/n}$ because $e\le n$.
Weyl's inequality therefore places one eigenvalue of $A_t/n$ near
$(k-1)/k$ and $k-1$ eigenvalues near $-1/k$, uniformly in all allowed
choices. Increase the threshold so that all these $k$ eigenvalues have
modulus at least $1/(2k)$; every other eigenvalue of $A_t$ has modulus at
most $\|tE\|\le\sqrt{2e}$. The two groups are disjoint.

If $v$ is a normalized eigenvector in the large group, then
\[
 \lambda(I-P_0)v=t(I-P_0)Ev.
\]
Summing the squared estimates over an orthonormal basis of that group
and using the equal ranks of $P_t$ and $P_0$ gives
\[
 \|P_t-P_0\|_F
 =\sqrt2\|(I-P_0)P_t\|_F
 \le\frac{2\sqrt2\,k}{n}\|EP_t\|_F.
\]
Consequently
\[
 \|EP_t\|_F\le\|EP_0\|_F+
 \frac{2\sqrt2\,k}{n}\|E\|\|EP_t\|_F
 \le\|EP_0\|_F+\frac{4k}{\sqrt n}\|EP_t\|_F.
\]
After additionally requiring $n\ge64k^2$, the last term is absorbed
and~\eqref{eq:localized-E} follows. This numerical condition
concerns only the absorption step; the earlier spectral-separation
threshold is still required.

Choose disjoint fixed neighborhoods of $(k-1)/k$, $-1/k$ and $0$ containing
the corresponding spectra of $A_t/n$. Let $\varphi_p$ be smooth, equal to
$|x|^p$ on the first two neighborhoods and zero on the third, with cutoffs
independent of $p$. Its first two derivatives are uniformly bounded on a
fixed compact interval for $p\in J$. Set
\[
 f_{n,p}(x)=n^p\varphi_p(x/n),\qquad
 F(t)=\tr f_{n,p}(A_t)=S_{p,k}(A_t)^p.
\]
In an orthonormal eigenbasis of $A_t$, the trace derivative formula
\cite{LewisSendov2001} reads
\[
 F'(t)=\tr(f'_{n,p}(A_t)E),\qquad
 F''(t)=\sum_{i,j}f_{n,p}'{}^{[1]}(\lambda_i,\lambda_j)|E_{ij}|^2.
\]
At repeated arguments, the divided difference is
$f'{}^{[1]}(x,x)=f''(x)$, so the formula also covers multiple eigenvalues.
The coefficient vanishes when both indices belong to the small group;
otherwise its absolute value is at most $C_{k,J}n^{p-2}$. Let $\mathcal L$
be the index set of the large group in this eigenbasis. Since $E$ is
symmetric,
\[
 \sum_{\substack{i,j\\i\in\mathcal L\ \mathrm{or}\ j\in\mathcal L}}
 |E_{ij}|^2
 =2\|EP_t\|_F^2-\|P_tEP_t\|_F^2
 \le2\|EP_t\|_F^2.
\]
Consequently,
\[
 |F''(t)|\le C_{k,J}n^{p-2}\|EP_t\|_F^2,\qquad
 |F'(t)|\le C_{k,J}n^{p-1}\|EP_t\|_F,\qquad
 F(t)\asymp_{k,J}n^p.
\]
The chain rule gives
\[
 \frac{d^2}{dt^2}F(t)^{1/p}
 =\frac1pF(t)^{1/p-1}F''(t)
 +\frac1p\left(\frac1p-1\right)F(t)^{1/p-2}F'(t)^2.
\]
The preceding bounds make both terms at most
$C_{k,J}n^{-1}\|EP_t\|_F^2$ in absolute value, proving
\eqref{eq:localized-second}. For $p=1$ the second term is zero;
the large eigenvalues remain separated from zero throughout.
\end{proof}

\subsection{The graph-theoretic loss formula}

\begin{lemma}\label{lem:sparse-loss}
Fix $k\ge4$ and the interval $J$ above.
Let $H=K_{n_1,\ldots,n_k}$ and let $G$ be obtained by adding $e_+$ edges within its parts and deleting $e_-$ edges between its parts.
Put
\[
 V=\sum_i(n_i-n/k)^2,\qquad e=e_++e_-.
\]
For the edit graph, write
\[
 d_E(v)=|\{w:vw\in E(G)\mathbin{\triangle}E(H)\}|,
 \qquad d_E=\max_v d_E(v).
\]
Uniformly for $p\in J$, $V+e\le n$ and sufficiently large $n$,
\begin{equation}\label{eq:sparse-loss}
 S_{p,k}(G)=nc_T(p)
 -\frac{g_k(p)V+z_+(p)e_++z_-(p)e_-}{n}
 +O_{k,J}\left(\frac{V^{3/2}+e\sqrt V+e d_E}{n^2}\right).
\end{equation}
In particular, $d_E=0$ when $e=0$.
\end{lemma}

\begin{proof}
The contribution of the part sizes is \eqref{eq:general-part-loss}.
At $A_0=A(H)$ the gradient on the selected spectral group is
\[
 X_0=U\mathcal G_p(Q(x))U^T,\qquad
 \mathcal G_p(Q)=S_{p,k}(Q)^{1-p}\sgn(Q)|Q|^{p-1}.
\]
At $x=k^{-1}\one$,
\[
 \mathcal G_p(Q)=D_p^{1-p}
 \left(\frac{(k-1)^{p-1}+1}{k}J_k-I_k\right).
\]
The function and the normalizing factors $x_i^{-1/2}$ are smooth on the fixed neighborhood.
Thus each lifted entry changes by $O_{k,J}(\sqrt V/n^2)$ when $x_i=1/k+\delta_i/n$.
For distinct vertices $v,w$,
\[
 (X_0)_{vw}=
 \begin{cases}
 -z_+(p)/(2n)+O_{k,J}(\sqrt V/n^2),&v,w\text{ in the same part},\\
 z_-(p)/(2n)+O_{k,J}(\sqrt V/n^2),&v,w\text{ in different parts}.
 \end{cases}
\]
Every added or deleted edge contributes twice the corresponding entry, with the appropriate sign.
Taylor's formula and Proposition~\ref{prop:localized} give
\[
 S_{p,k}(G)-S_{p,k}(H)
 =-\frac{z_+(p)e_++z_-(p)e_-}{n}
 +O_{k,J}\bigl((e\sqrt V+ed_E)/n^2\bigr).
\]
Combining with \eqref{eq:general-part-loss} proves the result.
\end{proof}
\begin{corollary}\label{thm:local-loss}
Fix $k\ge4$ and $B_1$.
If $|n_i-n/k|\le B_1$ and $G$ differs from $H=K_{n_1,\ldots,n_k}$ in at most $B_1$ edges, then, writing $e=|E(G)\triangle E(H)|$,
\begin{equation}\label{eq:local-loss}
 S_{\vartheta_k,k}(G)=nC_k-\frac{\gamma_k\sum_i(n_i-n/k)^2+\zeta_ke}{n}
             +O_{k,B_1}(n^{-2}).
\end{equation}
Every sequence with deficit $O(n^{-1})$ has such a representation with uniform bounds.
\end{corollary}
\begin{proof}
Use Lemma~\ref{lem:sparse-loss} at $\vartheta_k$, where all deviation and edit quantities are bounded.
The representation follows from Corollary~\ref{thm:bounded-edit}.
\end{proof}

\begin{corollary}\label{cor:vanishing}
Fix $k\ge4$ and let $G_n$ have order $n\to\infty$.
Then
\[
 \Delta_n:=nC_k-S_{\vartheta_k,k}(G_n)=o(1)
\]
if and only if there are complete $k$-partite graphs $H_n=K_{n_1,\ldots,n_k}$ on the same vertex sets such that
\[
 e_n:=|E(G_n)\mathbin{\triangle}E(H_n)|=o(n),\qquad
 V_n:=\sum_i(n_i-n/k)^2=o(n).
\]
For the representation supplied by Theorem~\ref{thm:multiscale},
\begin{equation}\label{eq:vanishing-loss}
 \Delta_n=\frac{\gamma_k V_n+\zeta_ke_n}{n}
 +O_k\left(\frac{\Delta_n^{3/2}}{\sqrt n}+\Delta_n^2\right).
\end{equation}
In particular, when $\Delta_n>0$ tends to zero, the quotient of $\Delta_n$ by $(\gamma_kV_n+\zeta_ke_n)/n$ tends to one.
For every fixed $0<\alpha\le1$, $\Delta_n=O(n^{-\alpha})$ implies
\[
 e_n=O_k(n^{1-\alpha}),\qquad
 \max_i|n_i-n/k|=O_k(n^{(1-\alpha)/2}),
\]
where the constants also depend on the constant in the assumed deficit bound.
\end{corollary}
\begin{proof}
If $\Delta_n=o(1)$, Theorem~\ref{thm:multiscale} gives $V_n+e_n\le K n\Delta_n=o(n)$.
Apply Lemma~\ref{lem:sparse-loss} at $p=\vartheta_k$ and use $d_E\le e_n$.
Its remainder is at most
\[
 C_k\frac{V_n^{3/2}+e_n\sqrt{V_n}+e_n^2}{n^2}
 \le C_k\left(\frac{\Delta_n^{3/2}}{\sqrt n}+\Delta_n^2\right),
\]
proving \eqref{eq:vanishing-loss}.
If $V_n+e_n=o(n)$ is given instead, the same lemma gives
\[
 \Delta_n=\frac{\gamma_kV_n+\zeta_ke_n}{n}
 +o\bigl((V_n+e_n)/n\bigr)=o(1).
\]
The zero case is interpreted by the exact formula when $V_n=e_n=0$.
\end{proof}

\subsection{Absolute deficits and failure of relative stability}\label{sec:obstruction}

\begin{theorem}\label{thm:nonstability}
Fix $k\ge4$ and $p=\vartheta_k$.
For every positive integer $t$, there is a connected graph $G_t$ of order $n=kt$ such that
\begin{equation}\label{eq:near-optimal}
 0\le nC_{p,k}-S_{p,k}(G_t)\le k^{1/p}.
\end{equation}
For every complete $k$-partite graph $H$ on $kt$ vertices, under every relabeling,
\[
\dist(G_t,H)\ge\frac12\left(\frac{t^2}{k}-t\right).
\]
In particular, $S_{p,k}(G_t)/(kt)\to C_{p,k}$, whereas
\[
 \liminf_{t\to\infty}\frac{\dist(G_t,T_k(kt))}{(kt)^2}\ge\frac1{2k^3}>0.
\]
\end{theorem}
\begin{proof}
Let $P$ be the symmetric permutation matrix interchanging the first two coordinates and fixing the others, and put $B=J_k-P$.
On the all-ones vector, $B$ has eigenvalue $k-1$.
On its orthogonal complement, $B=-P$ and all eigenvalues have modulus one.
Therefore
\[
 S_{p,k}(B)^p=(k-1)^p+k-1.
\]
Let $D_0$ be the diagonal matrix with the same diagonal as $B$, and define
\[
 A_t=B\otimes J_t-D_0\otimes I_t.
\]
This is a symmetric $(0,1)$ matrix with zero diagonal.
Its graph $G_t$ has $k$ classes of size $t$: the first two induce cliques and have no edges between them; the remaining $k-2$ classes are independent; all other pairs of classes are completely joined.
It is connected because $k\ge4$.

The nonzero singular values of $B\otimes J_t$ are $t$ times those of $B$.
Hence, at $p=\vartheta_k$,
\[
 S_{p,k}(B\otimes J_t)=t\binom{k}{2}^{1/p}=nC_{p,k}.
\]
The reverse triangle inequality and $S_{p,k}(D_0\otimes I_t)\le k^{1/p}$ yield the upper bound on the norm difference in \eqref{eq:near-optimal}.
Its nonnegativity follows from Theorem~\ref{thm:critical}.

Let $U$ be the first clique of $G_t$, and let the independent classes of $H$ meet $U$ in sizes $a_1,\ldots,a_k$.
Every edge within one of these intersections must be deleted.
Since $\sum_i a_i=t$, Cauchy's inequality gives
\[
 \dist(G_t,H)\ge\sum_{i=1}^k\binom{a_i}{2}
 =\frac12\left(\sum_i a_i^2-t\right)
 \ge\frac12\left(\frac{t^2}{k}-t\right).
\]
This bound holds for every partition and therefore for every relabeling.
Divide by $(kt)^2$ and let $t\to\infty$.
\end{proof}

\begin{remark}\label{rem:four-class-exact}
For $k=4$ and $t\ge2$, the preceding construction has
\[
 S_{1,4}(G_t)=\sqrt{16t^2+1}+2t-1,\qquad
 6t-S_{1,4}(G_t)=1-\frac1{\sqrt{16t^2+1}+4t}\longrightarrow1.
\]
Indeed, its four-dimensional quotient has characteristic polynomial
$(z+t)(z-t+1)(z^2-(2t-1)z-3t^2-t)$; the remaining eigenvalues are
$-1$ and $0$. Nevertheless its distance from every complete four-partite
graph is at least $t^2/8-t/2$.
\end{remark}

\subsection{The maximum and the second-largest value}\label{sec:order-extrema}
In this subsection put $p=\vartheta_k$ and $C=C_k$.
Write
\[
 M_{p,k}(n)=\max_{|V(G)|=n}S_{p,k}(G),\qquad
 M^{(2)}_{\vartheta_k,k}(n)=\max_{\substack{|V(G)|=n\\G\not\cong T_k(n)}}S_{\vartheta_k,k}(G).
\]

We use the uniform remainders in the following order. A candidate with
deficit at most $K_0/n$ confines every competing maximizer, by
Corollary~\ref{thm:bounded-edit}, to a common bound on the part deviations
and edit count. Corollary~\ref{thm:local-loss} then gives a common
remainder bound $C/n^2$. A positive leading loss $a_0/n$ is decisive for
$n>C/a_0$, in addition to the preceding order thresholds. All these
constants are chosen before the competing graph.

\begin{theorem}\label{thm:eventual}
For every fixed integer $k\ge4$, there exists $N_k$ such that for every $n\ge N_k$,
\[
 S_{\vartheta_k,k}(G)\le S_{\vartheta_k,k}(T_k(n))
 \qquad (|V(G)|=n),
\]
with equality if and only if $G\cong T_k(n)$.
Consequently, if $n\equiv r\pmod k$, $0\le r<k$, then
\[
 M_{\vartheta_k,k}(n)=nC_k-\frac{\gamma_k r(k-r)}{kn}+O_k(n^{-2}).
\]
\end{theorem}
\begin{proof}
Let $G$ be a maximizer on $n$ vertices.
Since $S_{p,k}(G)\ge S_{p,k}(T_k(n))$, Lemma~\ref{lem:secondorder} gives its norm deficit as $O_k(n^{-1})$.
Corollary~\ref{thm:bounded-edit} supplies a complete $k$-partite graph $H=K_{n_1,\ldots,n_k}$ with uniformly bounded part deviations and at most a bounded number $e$ of edits.
Let
\[
 V=\sum_i(n_i-n/k)^2,\qquad v_r=r(k-r)/k.
\]
Among integer $k$-tuples of sum $n$, $V$ is minimized precisely by the balanced tuple, with value $v_r$.
Moreover, $V-v_r$ is a nonnegative even integer: the corresponding sums of squares are integers congruent to $n$ modulo two.
Corollary~\ref{thm:local-loss} and Lemma~\ref{lem:secondorder} therefore give
\[
 0\le S_{p,k}(G)-S_{p,k}(T_k(n))
 =-\frac{\gamma_k(V-v_r)+\zeta_k e}{n}+O_k(n^{-2}).
\]
If $e\ge1$ or $V>v_r$, the numerator is at least $\min\{\zeta_k,2\gamma_k\}>0$.
Uniformity of the remainder makes this impossible for all sufficiently large $n$.
Thus $e=0$ and $V=v_r$, so $G\cong T_k(n)$.
The asymptotic formula follows from Lemma~\ref{lem:secondorder}.
\end{proof}

\begin{theorem}\label{thm:runnerup}\label{cor:second-largest}
Fix $k\ge4$. For all sufficiently large $n=kt+r$, the second-place graph is unique up to isomorphism. Define the available part-size types
\[
 \mathcal A_{k,r}(t)=\{t-1,t^{[k-r-2]},(t+1)^{[r+1]}\}\quad(k-r\ge2),
\]
\[
 \mathcal B_{k,r}(t)=\{t^{[k-r+1]},(t+1)^{[r-2]},t+2\}\quad(r\ge2).
\]
For $r=0,1$ the unique second-place graph has type $\mathcal A_{k,r}(t)$; for $r=k-1$ it has type $\mathcal B_{k,r}(t)$. When $2\le r\le k-2$, it has type $\mathcal B_{k,r}(t)$ for $k=4,5$, and type $\mathcal A_{k,r}(t)$ for $k\ge6$. Writing $M^{(2)}(n)=M^{(2)}_{\vartheta_k,k}(n)$, the gap is
\[
 M_{\vartheta_k,k}(n)-M^{(2)}(n)=
 \begin{cases}
 \displaystyle \frac{2\gamma_k}{n}+\frac{6\xi_k r}{kn^2}+O_k(n^{-3}),
       &\text{for type }\mathcal A_{k,r}(t),\\[5pt]
 \displaystyle \frac{2\gamma_k}{n}-\frac{6\xi_k(k-r)}{kn^2}+O_k(n^{-3}),
       &\text{for type }\mathcal B_{k,r}(t).
 \end{cases}
\]
\end{theorem}

\begin{proof}

Each available displayed complete multipartite graph has $V=v_r+2$, so they give a second-place candidate with reciprocal-order deficit.
Any second-place maximizer consequently satisfies Corollary~\ref{thm:bounded-edit} and \eqref{eq:local-loss}.
The constants satisfy
\[
 \zeta_k-2\gamma_k=(kC_k)^{1-\vartheta_k}[1+(\vartheta_k-1)(k-2)]>0.
\]
If its representation has $e\ge1$, then its loss coefficient is at least $\gamma_kv_r+\zeta_k$, strictly larger than $\gamma_k(v_r+2)$.
If $e=0$ but $V\ge v_r+4$, its coefficient is again strictly larger.
For large $n$, both possibilities are excluded by the uniform remainder.
The case $e=0,V=v_r$ is the excluded graph $T_k(n)$.
Thus $e=0,V=v_r+2$.
Writing $n_i=t+z_i$, this condition becomes
\[
 \sum_i z_i=r,\qquad\sum_i z_i(z_i-1)=2.
\]
Every summand is a nonnegative even integer, so precisely one $z_i$ equals $-1$ or $2$, and all others are zero or one.
These are exactly the two displayed multisets.
At least one type is available for every residue.

The preceding argument restricts every sufficiently large second-place graph to the displayed types. Relative to $T_k(n)$, both increase $\sum_i\delta_i^2$ by $2$. Their changes in $\sum_i\delta_i^3$ are $-6r/k$ and $6(k-r)/k$, respectively: write $n_i=t+z_i$, expand $(z_i-r/k)^3$, and use the displayed multiplicities. Therefore, when both types exist,
\[
 S_{\vartheta_k,k}(K_{\mathcal B_{k,r}(t)})-
 S_{\vartheta_k,k}(K_{\mathcal A_{k,r}(t)})
 =\frac{6\xi_k}{n^2}+O_k(n^{-3}).
\]
Here $K_{\mathcal A}$ denotes the complete multipartite graph with the specified part sizes. The sign assertion of Lemma~\ref{thm:cubic} gives the strict ordering and excludes an eventual tie. The single available type wins at the remaining residues. Subtracting \eqref{eq:cubic-expansion} from the balanced expansion proves the gap formulas.
\end{proof}

\begin{remark}\label{rem:small-order-runnerup}
The order restriction in Theorem~\ref{thm:runnerup} cannot be omitted.
For $k=4$ and $n=7$, the stated part-size rule gives
$H=K_{3,2,1,1}$. Let $G=\overline{P_3\bigsq K_2\bigsq K_2}$.
Both graphs have $17$ edges and are different from $T_4(7)$. Their
characteristic polynomials are
\[
 \begin{aligned}
 \chi_G(x)&=x^2(x+1)(x+2)(x^3-3x^2-10x+4),\\
 \chi_H(x)&=x^3(x+1)(x^3-x^2-16x-18).
 \end{aligned}
\]
Evaluation at rational endpoints places the roots of the first cubic in
$(-23/10,-11/5)$, $(0,1)$ and $(24/5,5)$, and those of the second in
$(-3,-2)$, $(-2,-1)$ and $(4,5)$. Thus the four largest singular values of
$G$ have sum greater than $24/5+11/5+2+1=10$. The graph $H$ has rank
four and only one positive eigenvalue, which is less than $5$; zero trace
therefore gives $S_{1,4}(H)<10$. Consequently,
\[
 S_{1,4}(G)>10>S_{1,4}(K_{3,2,1,1}).
\]
This comparison excludes $H$ as the second-place graph at order seven;
it does not identify all graphs attaining the second-largest value at
that order.
\end{remark}

\begin{corollary}\label{cor:xi-four}
For all sufficiently large $n=4t+r$, $0\le r\le3$,
\[
 M_{1,4}(n)=
 \begin{cases}
 6t,&r=0,\\
 2t+r-1+\sqrt{(2t+r-1)^2+12t(t+1)},&1\le r\le3.
 \end{cases}
\]
The unique maximizing graph is $T_4(n)$. For $r=0$ the value $6t$
holds for every $t\ge1$, as already recorded in~\cite[p.~2399]{Nik2011}.
\end{corollary}
\begin{proof}
By Theorem~\ref{thm:eventual}, it suffices to evaluate $S_{1,4}(T_4(n))$.
A complete four-partite graph with nonempty parts has rank four and one
positive eigenvalue, hence its norm is twice its spectral radius. If
$r>0$, grouping the $r$ parts of size $t+1$ and the $4-r$ parts of size
$t$ gives the equitable quotient
\[
 Q=\begin{pmatrix}
 (r-1)(t+1)&(4-r)t\\
 r(t+1)&(3-r)t
 \end{pmatrix}.
\]
Its Perron root is that of $T_4(n)$, since the Perron vector is constant
on each of these two vertex classes. The trace is $2t+r-1$ and the
determinant is $-3t(t+1)$. Solving the resulting quadratic gives the
formula. For $r=0$, the graph is $3t$-regular and its norm is $6t$.
\end{proof}

\section{Extremal graphs with a fixed edge deficit or surplus}\label{sec:prescribed}
\subsection{A finite family of part-size deviations}
Write $e_k(n)=|E(T_k(n))|$. Fix $k\ge4$, $0\le r<k$ and an integer $\ell\ge0$,
and define the finite set of integer multisets
\[
 \mathcal Z_{k,r,\ell}=
 \left\{(z_1,\ldots,z_k)\in\mathbb Z^k:
       \sum_i z_i=r,\quad\sum_i z_i(z_i-1)=2\ell\right\}/\mathfrak S_k.
\]
Each summand $z_i(z_i-1)$ is nonnegative. Thus the set is finite and independent
of $n$. For $n=kt+r$ and all sufficiently large $t$, its associated graph family is
\[
 \mathcal H_{k,r,\ell}(t)=
 \{K_{t+z_1,\ldots,t+z_k}:z\in\mathcal Z_{k,r,\ell}\}.
\]
Every member has exactly $e_k(n)-\ell$ edges. Balancing the parts maximizes the
number of edges \cite[Section~7.1]{Diestel2000}; here we instead prescribe the
integer loss in that number, so ordinary balancing would leave the admissible
class.

\begin{theorem}
\label{thm:deficit-shell}
Suppose $\mathcal Z_{k,r,\ell}\ne\varnothing$. For all sufficiently large
$n=kt+r$, every member of $\mathcal H_{k,r,\ell}(t)$ has strictly larger
$S_{\vartheta_k,k}$ than every $n$-vertex graph with $e_k(n)-\ell$ edges outside this
family. In particular, every maximizer is complete $k$-partite and belongs to
$\mathcal H_{k,r,\ell}(t)$.

Among the multisets in $\mathcal Z_{k,r,\ell}$, let $\mathcal Z^*$ maximize
$\sum_i z_i^3$ when $k=4,5$, and minimize it when $k\ge6$. Every maximizer
belongs to the subfamily indexed by $\mathcal Z^*$. If $\mathcal Z^*$ consists
of one multiset, the maximizing graph is unique up to isomorphism.
\end{theorem}
\begin{proof}
Put $v_r=r(k-r)/k$. If $n_i=t+z_i$, a direct calculation gives
\[
 V:=\sum_i(n_i-n/k)^2=v_r+2\ell,
 \qquad |E(K_{n_1,\ldots,n_k})|=e_k(n)-\ell.
\]
All graphs in the finite family have bounded part deviations. By
Corollary~\ref{thm:local-loss}, their deficits from $nC_k$ are
$O_{k,r,\ell}(n^{-1})$. Fix a common deficit bound for this finite family;
the reconstruction and remainder constants below then depend only on
$k,r,\ell$, not on its chosen member or the competing graph. Suppose that a graph $G$ outside this family were at least as good as one of its
members. The same upper bound on its deficit holds, uniformly in the competing
member. Corollary~\ref{thm:bounded-edit} supplies a complete $k$-partite graph $H$
with bounded part deviations, from which $G$ is obtained by $e_+$ internal
additions and $e_-$ cross-part deletions, with $e=e_++e_-$ bounded.

Let $V$ now denote the part variance of $H$, and put $w=(V-v_r)/2$.
The integer parity and balancing argument gives $w\in\mathbb Z_{\ge0}$.
The edge budget is the exact identity
\[
 e_- - e_+=\ell-w,\qquad e\ge|\ell-w|.
\]
Relative to the common leading loss coefficient of the shell, the additional
loss in Corollary~\ref{thm:local-loss} is
\[
 2\gamma_k(w-\ell)+\zeta_k e.
\]
If $w<\ell$, it is at least $(\zeta_k-2\gamma_k)(\ell-w)>0$; if $w>\ell$,
it is at least $(\zeta_k+2\gamma_k)(w-\ell)>0$. We used
\[
 \zeta_k-2\gamma_k=(kC_k)^{1-\vartheta_k}
       [1+(\vartheta_k-1)(k-2)]>0.
\]
If $w=\ell$ and $e>0$, it is again positive, at least $\zeta_k$.
These are fixed positive gaps between integers, whereas the norm remainder
is uniformly $O(n^{-2})$. Hence all three possibilities are impossible for
sufficiently large $n$. We must have $w=\ell$ and $e=0$, which places $G$ in
the shell, a contradiction. This proves the stronger first assertion, not just
the claim about a single maximizer.

Within the shell, Lemma~\ref{thm:cubic} applies uniformly. Since
$\sum z_i=r$ and $\sum z_i^2=r+2\ell$ are fixed,
\[
 \sum_i(z_i-r/k)^3
 =\sum_i z_i^3-\frac{3r}{k}(r+2\ell)+\frac{2r^3}{k^2}.
\]
Thus a strict integer gap in $\sum z_i^3$ gives a strict gap of order $n^{-2}$,
which dominates the uniform $O(n^{-3})$ remainder. The sign of $\xi_k$ in
Lemma~\ref{thm:cubic} proves the second assertion. Distinct multisets of
positive part sizes determine nonisomorphic complete multipartite graphs.
\end{proof}

\subsection{A family of unique extremizers}
\begin{lemma}\label{lem:integer-cubic-shell}
Let $1\le a\le k-1$ and $z_1,\ldots,z_k\in\mathbb Z$ satisfy
$\sum_i z_i=0$ and $\sum_i z_i^2=a(a+1)$. Then
\[
 -a(a^2-1)\le\sum_i z_i^3\le a(a^2-1).
\]
Equality on the left holds precisely for the multiset
$\{-a,1^{[a]},0^{[k-a-1]}\}$; equality on the right holds precisely for
$\{a,(-1)^{[a]},0^{[k-a-1]}\}$.
\end{lemma}
\begin{proof}
Every $|z_i|\le a$, because $(a+1)^2>a(a+1)$. For an integer $z\ge-a$,
\[
 z(z-1)(z+a)\ge0.
\]
Summing its expansion gives
$\sum z_i^3\ge-(a-1)\sum z_i^2=-a(a^2-1)$.
Equality requires $z_i\in\{-a,0,1\}$ for all $i$. If $b$ entries equal $-a$,
the zero sum forces $ab$ entries to equal $1$, and the second moment gives
$b(a^2+a)=a^2+a$. Thus $b=1$. Reflection proves the upper bound and its
unique equality multiset. This also covers $a=1$, when the two multisets agree.
\end{proof}

\begin{theorem}
\label{thm:triangular-deficit}
Fix $k\ge4$ and $1\le a\le k-1$, and put $\ell=a(a+1)/2$.
For all sufficiently large $n=kt$, the unique maximizer of $S_{\vartheta_k,k}$ over
all $n$-vertex graphs with $e_k(n)-\ell$ edges is
\[
 \begin{cases}
 K_{t+a,(t-1)^{[a]},t^{[k-a-1]}},&k=4,5,\\
 K_{t-a,(t+1)^{[a]},t^{[k-a-1]}},&k\ge6.
 \end{cases}
\]
The actual maximum equals
\begin{equation}\label{eq:triangular-deficit-max}
 nC_k-\frac{\gamma_k a(a+1)}n
     +\frac{|\xi_k|a(a^2-1)}{n^2}+O_{k,a}(n^{-3}).
\end{equation}
Repetitions of length zero are omitted from the displayed part sizes.
\end{theorem}
\begin{proof}
The displayed multisets show that the shell is nonempty. For $r=0$ its
constraints are exactly those of Lemma~\ref{lem:integer-cubic-shell}.
Theorem~\ref{thm:deficit-shell} and the sign of $\xi_k$ give uniqueness.
The variance and cubic sum of the winning deviations are $a(a+1)$ and
$\sgn(\xi_k)a(a^2-1)$, respectively. Substitution in
Lemma~\ref{thm:cubic} gives \eqref{eq:triangular-deficit-max}.
\end{proof}

\subsection{A fixed edge surplus}
A surplus requires actual internal edges after reconstruction.
Their degree-square sum governs the next loss term, while a mixed part-size term identifies the correct part at nondivisible orders.
Write $\kappa=k-1$, $D_p=(\kappa^p+\kappa)^{1/p}$, and set
\begin{align*}
 \beta_{p,k}&=kD_p^{1-p}(\kappa+\kappa^{p-2}),\\
 \chi_{p,k}&=2kD_p^{1-p}
 \left[k-(k-2)(p-1)-\frac{2(1+\kappa^{p-1})}{k}\right].
\end{align*}
The bracket in $\chi_{p,k}$ is zero at $p=2$ and has strictly negative derivative.
Thus both coefficients are positive on $[1,2)$, and
\begin{equation}\label{eq:chi-critical}
 \chi_k:=\chi_{\vartheta_k,k}=2k(kC_k)^{1-\vartheta_k}[1+(k-2)(2-\vartheta_k)]>0.
\end{equation}

\begin{lemma}\label{lem:internal-walk}
Fix $k\ge4$, a nonnegative integer $\ell$, a constant $B_0$, and a compact interval $J\subset[1,2)$. Let $H=K_{n_1,\ldots,n_k}$, where $\delta_i=n_i-n/k$, $\sum_i\delta_i=0$, and $|\delta_i|\le B_0$. Add $\ell$ internal edges to obtain $H+F$, let $\ell_i$ be their number in part $i$, and put $W(F)=\sum_vd_F(v)^2$. Uniformly over these choices and $p\in J$,
\begin{equation}\label{eq:internal-walk}
 S_{p,k}(H+F)-S_{p,k}(H)
 =-\frac{z_+(p)\ell}{n}
 +\frac{\beta_{p,k}W(F)+\chi_{p,k}\sum_i\delta_i\ell_i}{n^2}
 +O_{k,J,\ell,B_0}(n^{-3}).
\end{equation}
If $F_s$ and $F_\triangle$ are a star $K_{1,3}$ and a triangle added to the same part of $H$, then
\begin{equation}\label{eq:star-triangle}
 S_{p,k}(H+F_s)-S_{p,k}(H+F_\triangle)
 =\frac{6k^2D_p^{1-p}(\kappa-\kappa^{p-3})}{n^3}+O_{k,J,B_0}(n^{-4}).
\end{equation}
The leading coefficient in \eqref{eq:star-triangle} is positive throughout $J$.
\end{lemma}
\begin{proof}
We use the argument principle on separated positive and negative spectral clusters. This retains multiplicities and avoids differentiating individually ordered eigenvalues at a repeated root; see~\cite[Section~1.2]{Kressner2005} for the cluster viewpoint.

Set $T=n/k$, which need not be an integer. In part $i$, let $E_i$ be the adjacency matrix of the added graph and define
\[
 h_i(Tz)=\one^T(TzI-E_i)^{-1}\one,\qquad
 w_i=\one^TE_i^2\one.
\]
At most $2\ell$ vertices are incident with an added edge, and $\|E_i\|\le\ell$.
For $j\ge1$, this also gives
$|\one^TE_i^j\one|\le2\ell\,\ell^j$ (with the zero case interpreted directly).
Choose fixed circles $\Gamma_+$ and $\Gamma_-$ about $\kappa$ and $-1$ of
radius less than $1/4$, together with slightly larger disjoint annular
neighborhoods that avoid zero. On these neighborhoods, the Neumann series
for $(TzI-E_i)^{-1}$ and its first derivative converge uniformly once
$T$ exceeds a constant depending only on $\ell$. The resolvent series therefore gives, on fixed complex contours separated from zero,
\[
 h_i(Tz)=z^{-1}+\frac{\delta_i}{Tz}
   +\frac{2\ell_i}{T^2z^2}+\frac{w_i}{T^3z^3}+O_\ell(T^{-4}).
\]
All walk coefficients after the constant term are bounded independently of
the part sizes. Let $Z$ be the $n$-by-$k$ matrix of unnormalized part
indicators, $\mathcal E=\diag(E_1,\ldots,E_k)$, and $K=J_k-I_k$.
Then $A(H+F)=\mathcal E+ZKZ^T$. For
$Tz\notin\bigcup_i\Spec(E_i)$, the matrix determinant lemma gives
\[
 \begin{aligned}
 \det\bigl(TzI_n-A(H+F)\bigr)
 &=\left(\prod_{i=1}^k\det(TzI_{n_i}-E_i)\right)\mathcal D_F(z),\\
 \mathcal D_F(z)&=\det\bigl(I_k-\diag(h_i(Tz))K\bigr).
 \end{aligned}
\]
Here $Z^T(TzI_n-\mathcal E)^{-1}Z=\diag(h_i(Tz))$, and
$\det(I-XY)=\det(I-YX)$ allows the displayed order of the factors.
Let $\mathcal D_H$ denote the same expression with $E_i=0$. For sufficiently
large $T$, the product of block determinants has no zeros on or inside
$\Gamma_+$ and $\Gamma_-$. Thus $\mathcal D_F$ preserves the zeros and
their multiplicities in these circles. The unperturbed $k$-by-$k$ quotient of $H/T$ converges uniformly to $K$; the perturbation due to $F/T$ has norm $O(T^{-1})$. Thus the two circles contain exactly the $k$ largest absolute eigenvalues, one positive and $k-1$ negative, for sufficiently large $n$.

Put $L(z)=K(I-K/z)^{-1}$. Its diagonal and off-diagonal entries are
\[
 l(z)=\frac{\kappa z}{(z-\kappa)(z+1)},\qquad
 b(z)=\frac{z^2}{(z-\kappa)(z+1)},
\]
and
\[
 l(z)^2-b(z)^2=-\frac{z^2(z+\kappa)}{(z-\kappa)(z+1)^2}.
\]
The expansion of $\log\det(I-X)$ about $X=0$, applied on neighborhoods of these circles, now gives
\begin{align*}
 \log\frac{\mathcal D_F(z)}{\mathcal D_H(z)}
 ={}&-\frac{2\ell \kappa}{T^2z(z-\kappa)(z+1)}
     -\frac{W(F)\kappa}{T^3z^2(z-\kappa)(z+1)}\\
 &+\frac{2(\sum_i\delta_i\ell_i)(z+\kappa)}{T^3z(z-\kappa)(z+1)^2}
     +O_{k,\ell,B_0}(T^{-4}).
\end{align*}
For clarity, the mixed term before simplification is
$-2T^{-3}z^{-3}\sum_{i,j}\delta_i\ell_jL_{ij}(z)^2$; using $\sum_i\delta_i=0$ gives the displayed expression. This is the term absent at divisible balanced orders.

Use $f(z)=z^p$ near $\Gamma_+$ and $f(z)=(-z)^p$ near $\Gamma_-$. The change in the sum of $f$ over the enclosed scaled eigenvalues equals
\[
 -\frac1{2\pi\mathrm i}\int_{\Gamma_+\cup\Gamma_-}
 f'(z)\log\frac{\mathcal D_F(z)}{\mathcal D_H(z)}\,dz.
\]
Integration by parts in the argument principle justifies this identity.
On the fixed annular neighborhoods, $\mathcal D_H$ converges uniformly to
$\det(I-K/z)$, which is bounded away from zero there. Moreover,
$\mathcal D_F/\mathcal D_H=1+O_{k,\ell,B_0}(T^{-2})$.
The logarithm therefore has its unique branch near zero. The Taylor
remainder is uniform on those neighborhoods, and Cauchy's estimate gives
the same order for its derivative on the circles. The derivatives of
$f$ are uniform for $p\in J$. The required residues are
\begin{align*}
 \frac1{2\pi\mathrm i}\int f'(z)\frac{\kappa}{z^j(z-\kappa)(z+1)}\,dz
 &=\frac pk\bigl[\kappa^{p-j}+(-1)^j\kappa\bigr],\\
 \frac1{2\pi\mathrm i}\int f'(z)\frac{z+\kappa}{z(z-\kappa)(z+1)^2}\,dz
 &=p\left[\frac{(p-1)(k-2)}k-1+\frac{2(1+\kappa^{p-1})}{k^2}\right].
\end{align*}
Here and below these integrals are over $\Gamma_+\cup\Gamma_-$. The second line follows from the simple pole at $\kappa$ and double pole at $-1$; the derivatives there are $f'(\kappa)=p\kappa^{p-1}$, $f'(-1)=-p$, and $f''(-1)=p(p-1)$.

The unedited scaled power sum is $D_p^p+O(T^{-2})$. Indeed its first-order term is permutation-invariant and linear in $(\delta_i)$, hence a multiple of $\sum_i\delta_i=0$. Taking the $p$th root and multiplying by $T$ therefore changes the first two displayed perturbation orders by the factor $D_p^{1-p}/p$; all quadratic root terms affect only $O(T^{-3})$ in the norm. Substituting $n=kT$ gives \eqref{eq:internal-walk} and the stated coefficients.

For a star and triangle in the same part, their first two nonconstant walk moments agree:
\[
 \one^TE_s\one=\one^TE_\triangle\one=6,\qquad
 \one^TE_s^2\one=\one^TE_\triangle^2\one=12.
\]
Their third moments are $18$ and $24$, respectively, by
$\one^TE^3\one=2\sum_{uv\in E}d(u)d(v)$. Hence the leading difference in the two logarithms, triangle minus star, is
$-6\kappa/[T^4z^3(z-\kappa)(z+1)]$. Every term involving part imbalance in this difference is $O(T^{-5})$, because the previous moments coincide. The first residue formula with $j=3$ and the root expansion prove \eqref{eq:star-triangle}. All functions and contours can be chosen uniformly on the compact interval $J$, including $p=1$.
\end{proof}

\begin{lemma}[\cite{Das2003}]\label{lem:incident-edges}
If $F$ is a connected graph with $n$ vertices and $\ell\ge1$ edges, then $W(F)\le\ell(\ell+1)$ with equality for $n > 3$ if and only if $G \cong K_3$ or $G \cong K_{1,n-1}$.
\end{lemma}

Write $n=kt+s$, $0\le s<k$, and put
\[
 v_s=s(k-s)/k,\quad w_s=s(k-s)(k-2s)/k^2,\quad
 d_s=\begin{cases}0,&s=0,\\1-s/k,&s>0.\end{cases}
\]

\begin{theorem}\label{thm:fixed-surplus}
Fix integers $k\ge4$ and $\ell\ge1$.
For every sufficiently large $n$, the unique maximizer of $S_{\vartheta_k,k}$ among all graphs with $n$ vertices and $e_k(n)+\ell$ edges is $T_k(n;\ell)$, obtained by adding a star $K_{1,\ell}$ inside a largest part of $T_k(n)$.
\end{theorem}
\begin{proof}
The stated candidate has deficit $O_{k,\ell}(n^{-1})$ by Lemma~\ref{lem:internal-walk} and the multipartite expansion. Hence every constrained maximum has a representation, from Corollary~\ref{thm:bounded-edit}, relative to $H=K_{n_1,\ldots,n_k}$ with bounded $V=\sum_i(n_i-n/k)^2$ and bounded edit count. Let $e_+$ and $e_-$ denote internal additions and cross-part deletions. Since the part sizes are integers summing to $n$,
\[
 V\ge v_s,\qquad V-v_s\in2\mathbb Z_{\ge0}.
\]
The prescribed edge count gives the exact relations
\begin{equation}\label{eq:surplus-balance}
 e_+-e_-=\ell+\frac{V-v_s}2,\qquad
 e_++e_-=\ell+\frac{V-v_s}2+2e_-.
\end{equation}
Substitute in Corollary~\ref{thm:local-loss}. Compared with the candidate, the extra leading loss is
\[
 \frac{(\gamma_k+\zeta_k/2)(V-v_s)+2\zeta_ke_-}{n}.
\]
Positivity and discreteness force $V=v_s$ and $e_-=0$ for all sufficiently large $n$. Thus $H=T_k(n)$ and all $\ell$ added edges are internal; this is a conclusion about the global constrained maximum, not an assumed graph class.

For internal additions the two second-order terms obey
\[
 W(F)\le\ell(\ell+1),\qquad \sum_i\delta_i\ell_i\le\ell d_s.
\]
Both are maximized simultaneously by a star in a largest part. Their coefficients $\beta_{\vartheta_k,k},\chi_k$ are positive. If either bound is strict, its gap is at least one: the first is an integer, and the second gap is the number of edges placed in smaller parts. The resulting loss of order $n^{-2}$ dominates the uniform remainder. Lemma~\ref{lem:incident-edges} leaves only the star, or a triangle when $\ell=3$, in a largest part. Equation~\eqref{eq:star-triangle} excludes the triangle. All remaining placements are isomorphic. This proves the stated uniqueness.
\end{proof}

\section{Vertex partitions and lower bounds}\label{sec:partitions}
For $k=0$, put $S_{p,0}=0$.
The full-energy pinching result in~\cite{Akbari2025} motivates the following
truncated statement. The inequality without a remainder is the
rank-budget form of classical pinching: if
$A_{\mathcal P}=\bigoplus_i A(G[V_i])$ and $\sum_i k_i=k$, then
\[
 \sum_i S_{p,k_i}(G[V_i])^p
 \le S_{p,k}(A_{\mathcal P})^p\le S_{p,k}(A)^p.
\]
The first inequality allocates $k$ singular values among the diagonal
blocks, and the second is pinching for a unitarily invariant norm.
The refinement below adds a quadratic term with a best universal
coefficient and gives the equality conditions. It concerns total absolute
$p$-energy, not all of the separate positive and negative $p$-energy
statements in~\cite{Akbari2025}. 
For a positive semidefinite matrix $B$, a maximizing $k$-dimensional spectral subspace contains the eigenspaces above its $k$th largest eigenvalue and the required dimension within the eigenspace at that value.
This convention includes zero eigenvalues and all ties. For a unit vector $u$, its spectral measure for $A$ assigns mass $\|P_\lambda u\|_2^2$ to each eigenvalue $\lambda$, where $P_\lambda$ is the corresponding spectral projection.

\begin{theorem}\label{thm:pinching}
Let $G$ be a nonempty graph, and let $V(G)=V_1\dot\cup\cdots\dot\cup V_s$ be a partition into nonempty sets.
Choose $0\le k_i\le|V_i|$ with $\sum_i k_i=k\ge1$.
For each $i$, let the columns of $U_i$ be orthonormal eigenvectors of $A(G[V_i])$ corresponding to its $k_i$ largest absolute eigenvalues, and put
\[
 \mathcal I=\sum_{i=1}^s\Frob{A[V\setminus V_i,V_i]U_i}^2.
\]
For $p\ge1$,
\begin{equation}\label{eq:budget}
 S_{p,k}(G)^p\ge\sum_{i=1}^s S_{p,k_i}(G[V_i])^p.
\end{equation}
For $1<p\le2$, this strengthens to
\begin{equation}\label{eq:quant-pinch}
 S_{p,k}(G)^p\ge\sum_{i=1}^s S_{p,k_i}(G[V_i])^p
        +\frac{p(p-1)}2\rho^{p-2}\mathcal I.
\end{equation}
The coefficient $p(p-1)/2$ is the largest possible universal coefficient in \eqref{eq:quant-pinch}, even for connected simple graphs with two parts and $k_1=k_2=1$.
For $p>1$, equality in \eqref{eq:budget} holds if and only if the embedded selected vectors are eigenvectors of $A(G)$ and their eigenvalue moduli exhaust its largest $k$ singular values, allowing ties.
For $p=1$, equality holds if and only if their span is a maximizing $k$-dimensional spectral subspace for $|A|$ and the spectral measure of each vector is supported in either $[0,\infty)$ or $(-\infty,0]$.
For $1<p<2$, equality in \eqref{eq:quant-pinch} has the same characterization as equality in \eqref{eq:budget}; in particular, $\mathcal I=0$ is necessary.
For $p=2$, equality in \eqref{eq:quant-pinch} holds if and only if the selected span is a maximizing $k$-dimensional spectral subspace for $A^2$, without requiring $\mathcal I=0$.
\end{theorem}
\begin{proof}
The selected vectors, extended by zeros outside their parts, form an orthonormal family $u_1,\ldots,u_k$.
By the maximum principle for the positive semidefinite matrix $|A|^p$,
\[
 S_{p,k}(G)^p\ge\sum_{j=1}^k u_j^T|A|^p u_j
 \ge\sum_{j=1}^k|u_j^TAu_j|^p.
\]
The second step is scalar Jensen applied to the spectral measure of each $u_j$.
This proves \eqref{eq:budget}.
For $1<p\le2$, the function
\[
 f(t)=|t|^p-\frac{p(p-1)}2\rho^{p-2}t^2
\]
is convex on $[-\rho,\rho]$.
Indeed, its second derivative is nonnegative away from zero, and the first derivative is continuous at zero.
Applying Jensen to $f$ gives, for each unit vector $u$,
\[
 u^T|A|^pu-|u^TAu|^p
 \ge\frac{p(p-1)}2\rho^{p-2}
     \big(\|Au\|_2^2-(u^TAu)^2\big).
\]
If $u$ is an embedded local eigenvector from $V_i$, the expression in parentheses is exactly
$\|A[V\setminus V_i,V_i](u|_{V_i})\|_2^2$.
Summing proves \eqref{eq:quant-pinch}.
For $p>1$, strict convexity of $|t|^p$ says that equality in scalar Jensen holds precisely when each $u_j$ is a global eigenvector.
Equality in the maximum principle additionally selects the largest $k$ moduli, proving the characterization for $p>1$ in \eqref{eq:budget}.
At $p=1$, equality in scalar Jensen for the absolute value is equivalent to a spectral measure supported in either $[0,\infty)$ or $(-\infty,0]$; the maximum principle gives the other stated condition.
For $1<p<2$, the function $f$ above is strictly convex: its derivative is strictly increasing throughout $(-\rho,\rho)$, including across zero. Equality in its Jensen inequality therefore forces each spectral measure to be a point mass, giving the asserted characterization for \eqref{eq:quant-pinch}.
At $p=2$, $f$ is identically zero, and the right-hand side of \eqref{eq:quant-pinch} is exactly $\sum_j u_j^TA^2u_j$. Only equality in the maximum principle remains.

For sharpness, let $G_a$ consist of two copies of $K_a$ joined by a perfect matching, where $a\ge4$.
Use the two cliques as parts and their constant unit vectors, so $\rho=a$ and $\mathcal I=2$.
Its adjacency spectrum is
\[
 a,\ a-2,\ 0^{[a-1]},\ (-2)^{[a-1]}.
\]
Therefore the quotient of the difference in \eqref{eq:budget} by $\rho^{p-2}\mathcal I$ equals
\[
 \frac{a^p+(a-2)^p-2(a-1)^p}{2a^{p-2}}
 \longrightarrow\frac{p(p-1)}2.
\]
No larger universal coefficient can hold. When $p=2$, the quotient equals $1$ exactly.
\end{proof}

\section{Further questions}\label{sec:questions}
The extremal results above lead to questions about attainment of the
order bound, the range of the saturation exponents, prescribed edge
counts and the scale of reconstruction.

\begin{problem}\label{prob:attainment}
For $1\le p<2$ and fixed $k\ge2$, determine $\mathfrak m_{p,k}$ and
characterize the pairs for which $\mathfrak m_{p,k}=C_{p,k}$.
Which pairs admit equality $M_{p,k}(n)=nC_{p,k}$ at a finite order?
\end{problem}
Proposition~\ref{prop:arithmetic} constrains the second question, not
the first. In particular, finite-order nonattainment for nonsquare $k$
at $p=1$ does not prove the strict asymptotic gap asked for infinitely
many $k$ in~\cite[Conjecture~4.7]{Nik2015}.

\begin{problem}\label{prob:interval}
For fixed $k\ge4$, determine the set
\[
 \mathcal I_k=\{p\in[1,2): T_k(n)\text{ is the unique maximizer of }
 S_{p,k}\text{ for all sufficiently large }n\}.
\]
Does $\mathcal I_k$ contain a relative interval about $\vartheta_k$?
What are the endpoints of its component containing $\vartheta_k$?
\end{problem}
Theorem~\ref{thm:eventual} proves $\vartheta_k\in\mathcal I_k$.
There is a simple upper obstruction. The function
\[
 c_T(p)=\frac1k\big((k-1)^p+k-1\big)^{1/p}
\]
is strictly decreasing and equals $C_k>1$ at $\vartheta_k$, while
$c_T(2)=\sqrt{(k-1)/k}<1$. Thus there is a unique
$\pi_k\in(\vartheta_k,2)$ with $c_T(\pi_k)=1$.
For $p>\pi_k$, the limits
$S_{p,k}(T_k(n))/n\to c_T(p)<1$ and
$S_{p,k}(K_n)/n\to1$ imply $p\notin\mathcal I_k$.
This comparison does not identify the actual endpoint.

\begin{problem}\label{prob:sizes}
Classify the maximizers when $\mathcal Z_{k,r,\ell}$ is empty or has
several multisets with the same extremal cubic sum. For
$\ell=\ell(n)\to\infty$, determine the range in which the star-in-a-largest-part
conclusion remains valid, subject to
$\ell\le\lceil n/k\rceil-1$.
\end{problem}
Growing $\ell$ requires uniform estimates beyond the fixed-budget
expansions. Ties in the cubic term may require higher-order comparisons.

\begin{problem}\label{prob:scale}
Let $\delta_k^*$ be the supremum of all $\delta>0$ for which constants
$K,N$ exist such that every $n$-vertex graph, $n\ge N$, with
$\Delta_k(G)\le\delta$ admits a complete $k$-partite graph $H$ satisfying
\[
 |E(G)\mathbin{\triangle}E(H)|+
 \sum_i(n_i-n/k)^2\le Kn\Delta_k(G).
\]
Determine $\delta_k^*$. More generally, describe the possible structures
of sequences with $S_{\vartheta_k,k}(G_n)/n\to C_k$.
\end{problem}
Theorems~\ref{thm:multiscale} and~\ref{thm:nonstability} give
\[
 0<\delta_k^*\le (kC_k)^{1-\vartheta_k}.
\]
For the upper bound, use the one-transposition construction from
Theorem~\ref{thm:nonstability}. On the class-indicator space its adjacency
matrix is $tB-D_0$, where $B=J_k-P$ and
$D_0=\diag(1,1,0,\ldots,0)$. Since
\[
 B|B|^{\vartheta_k-2}=\tfrac12J_k-P,
 \qquad \tr\big(B|B|^{\vartheta_k-2}D_0\big)=1,
\]
first-order differentiation of the separated nonzero spectral clusters gives
\[
 S_{\vartheta_k,k}(G_t)
 =tkC_k-(kC_k)^{1-\vartheta_k}+O_k(t^{-1}).
\]
The remaining eigenvalues are bounded and do not enter the first $k$
moduli for large $t$. Together with the quadratic edit-distance bound in
Theorem~\ref{thm:nonstability}, this proves the upper bound on $\delta_k^*$.
Whether it is sharp remains open.

\begin{problem}\label{prob:minima}
For fixed $p\ge1$ and $2\le k<n$, determine
\[
 m^{\mathrm{conn}}_{p,k}(n)
 =\min\{S_{p,k}(G): |V(G)|=n,\ G\text{ connected}\}
\]
and characterize all minimizing graphs.
\end{problem}
This is a truncated analogue of~\cite[Question~4.51]{Nik2016}.
Tang, Liu and Wang~\cite{Tang2026} proved full-spectrum star-minimality
for $1<p<2$. For fixed $k$, 
\[
 S_{p,k}(P_n)<2k^{1/p},\qquad
 S_{p,k}({K_{1,n-1}})=2^{1/p}\sqrt{n-1}\quad(k\ge2).
\]
Consequently, the star is not a minimizer whenever
$n-1>2^{2-2/p}k^{2/p}$. The first inequality uses
$\rho(P_n)<2$, and the second uses the two nonzero star eigenvalues.

Explicit thresholds, small-order exceptions and further ranks in the
ordering at $\vartheta_k$ also remain to be determined. For the
full-Schatten maximum problem when $p>2$, see Bo\v zin and
Lazarevi\'c~\cite{Bozin2025}.

\end{document}